\documentclass[final]{amsart}
\usepackage{amssymb,amsmath}

\usepackage{amssymb}

\newtheorem{Theorem}{Theorem}[section]
\newtheorem{lemma}{Lemma}[section]
\newtheorem{corollaire}{Corollary}[section]

\theoremstyle{remark}
\newtheorem{remark}{Remark}[section]

\newcommand{\R}{{\mathbb R}}

\newcommand{\argmin}{{\rm argmin}\kern 0.12em}

\newcommand{\Z}{\mathcal Z}

\newcommand{\xb}{\overline x}

\newcommand{\Xiu}{\mathcal X_1}
\newcommand{\Xid}{\mathcal X_2}

\begin{document}

\title{Asymptotic behavior of coupled dynamical systems with multiscale aspects
  }

\author{Hedy Attouch}

\author{Marc-Olivier Czarnecki}

\address{Institut de Mathématiques et Modélisation de Montpellier, UMR 5149 CNRS, Universit\'e Montpellier 2, place Eug\`ene Bataillon,
34095 Montpellier cedex 5, France}
\email{attouch@math.univ-montp2.fr, marco@math.univ-montp2.fr}

 \thanks{with the support of the French ANR 
   under grants ANR-05-BLAN-0248-01 and ANR-08-BLAN-0294-03.}
\maketitle

\vspace{0.3cm}

\paragraph{\textbf{Abstract}}  We study the asymptotic behaviour, as  time variable $t$ goes to $+\infty$, of  nonautonomous dynamical systems involving multiscale features.
 As a benchmark case, given $\mathcal H$ a general Hilbert space, $\Phi: \mathcal H \rightarrow \R \cup \{+\infty\}$ and $\Psi: \mathcal H \rightarrow \R \cup \{+ \infty\}$  two closed convex functions, and $\beta$  a function of $t$ which tends to $+\infty$ as $t$ goes to $+\infty$, 
we consider the  differential inclusion
$$\dot{x}(t) +  \partial \Phi (x(t)) +  \beta (t) \partial \Psi (x(t))\ni 0. $$
This system  models the emergence of various  collective behaviors in game theory, as well as the asymptotic control of coupled systems. 
We show several results ranging from weak ergodic to strong convergence of the trajectories. As a key ingredient we 
assume that, for every $p$ belonging to the range of $N_C$ 
$$  \int_{0}^{+\infty} \beta (t) \left[\Psi^* \left(\frac{p}{ \beta (t)}\right) -
\sigma_C \left(\frac{p}{ \beta (t)}\right)\right]dt < + \infty
$$
where $\Psi^*$ is the Fenchel conjugate of $\Psi$, $ \sigma_C $ is the support function of $C= \mbox{argmin} \Psi$ and
 $ N_C (x)$ is the normal cone to $C$ at $x$.
As a by-product, we revisit the  sytem
$$
\dot{x}(t) +  \epsilon (t) \partial \Phi (x(t)) +  \partial \Psi
(x(t))\ni 0 
$$
where  $\epsilon (t)$ tends to zero as $t$ goes to $+\infty$ and $  \int_{0}^{+\infty}\epsilon (t) dt = + \infty  $, whose asymptotic behaviour can be  derived from the preceding one by  time rescaling.
Applications are given in game theory, optimal control, variational problems and PDE's.

\vspace{0.2cm}


\vspace{0.2cm}


\markboth{H. ATTOUCH, M.-O. CZARNECKI}
  {MULTISCALED GRADIENT DYNAMICS}

\indent

\newpage

\section{Introduction}

$\mathcal H$ is a real Hilbert space, we write $ \|x\|^2 = \left\langle x , x\right\rangle$ for $x\in \mathcal H$.
We denote by $\Gamma_0(\mathcal H)$ the class of   closed (lower semicontinuous) convex proper functions 
from $\mathcal H$ to $\R \cup \{+\infty\}$ which are not identically equal to $+\infty$. The 
subdifferential of $f\in\Gamma_0(\mathcal H)$ is the maximal monotone operator

\begin{eqnarray*}
\label{e:subdiff}
\partial f\colon \mathcal H&\to& 2^{\mathcal H}\\
x&\mapsto&
\left\{u\in\mathcal H\ | \ ( \forall y\in\mathcal H)\;\; f(y) \geq  f(x) + \left\langle u, y - x \right\rangle  \right\}.\nonumber
\end{eqnarray*}

\subsection{Problem statement} 

\begin{itemize}
	
\item $\Phi: \mathcal H \rightarrow \R^{\phantom{+}} \cup \{+\infty\}$ is a closed convex proper function.

\item $\Psi: \mathcal H \rightarrow \R^+ \cup \{+ \infty\}$  is a closed convex proper function, $C$ = $\mbox{argmin}\Psi = \Psi ^{-1}(0) \neq \emptyset$.
	
\item  $\beta : \R^+ \rightarrow \R^+$ is  a function of $t$ which tends to $+\infty$ as $t$ goes to $+\infty$.

\end{itemize}

We study the asymptotic behavior of the trajectories of the nonautonomous multiscaled differential inclusion
\begin{equation}\label{MAG1}
(MAG) \qquad \dot{x}(t) +  \partial \Phi (x(t)) +  \beta (t) \partial \Psi (x(t))\ni 0
\end{equation}

where $\partial \Phi$ and $\partial \Psi$ are the subdifferentials  of $\Phi$ and $\Psi$. 

\vspace{0.1cm}

Let us observe that $\partial \Phi +  \beta (t) \partial \Psi \subset \partial (\Phi + \beta (t) \Psi)$   (equality holds under some general qualification assumption). Hence, each trajectory of $(MAG)$ satisfies  

$$ \qquad \dot{x}(t) +   \partial \left(\Phi + \beta (t)  \Psi  \right)(x(t))\ni 0.$$

On the other hand, $\mbox{ }$ $ \Phi + \beta(t)\Psi \uparrow \Phi + \delta_C$ as $t\rightarrow +\infty$
where $\delta_C$ is the indicator function of the set $C$ ($\delta_C (x) = 0$  for $x\in C, +\infty$ outwards). 
Monotone convergence is a variational convergence, (\cite{A}, theorem 3.20). As a consequence, 
 the corresponding subdifferential operators converge in the sense of graphs (equivalently in the sense of resolvents) as $t\rightarrow +\infty$ (\cite{A}, theorem 3.66)

$$ \partial \left(\Phi + \beta (t)  \Psi  \right) \rightarrow \partial (\Phi + \delta_C).$$

From the asymptotical point of view, this suggests strong  analogies between  $(MAG)$ and the steepest descent dynamical system associated to the closed convex proper function
$\Phi + \delta_C \in \Gamma_0(\mathcal H)$

\begin{equation}\label{MAG2}
\qquad \dot{x}(t) +  \partial \left(\Phi + \delta_C \right)(x(t))\ni 0.
\end{equation}

In our main result, theorem 3.1, we  prove that the two systems (\ref{MAG1}) and (\ref{MAG2}) share  similar asymptotical properties, whence the terminology $(MAG)$= ``Multiscale Asymptotic Gradient'' system. 

More precisely, under general assumptions,
we  prove that  each trajectory of  $(MAG)$  weakly converges in $\mathcal H$, with its limit belonging to $\mbox{argmin}_C \Phi$

\begin{equation}\label{MAG3}
x(t)\rightharpoonup x_{\infty} \in \mbox{argmin}_C \Phi  \mbox{ as }  t\rightarrow +\infty.
\end{equation}

This result can be seen as an extension of Bruck's theorem \cite{Bruck} to multiscaled nonautonomous gradient systems. 

\subsection{Notion of solution.}

We  consider strong solutions in the sense of Brezis (\cite{Bre1}, definition 3.1). 
Such a solution $x(.)$ is continuous on $\left[0, + \infty \right)$ and absolutely continuous on any bounded interval
 $\left[0,T\right]$ with $T< +\infty.$
 Being absolutely continuous, $x(.)$ is almost everywhere
 differentiable, and it it assumed that the equation holds almost
 everywhere. Equivalently $x(.)$ is a solution of $(MAG)$ if
there exist two functions $\xi(.)$ and $\eta (.)$ with

\begin{equation*}\label{strongsol02}
\xi(t) \in \partial \Phi(x(t)) \mbox{ }\mbox{and}\mbox{ } \eta (t) \in \partial \Psi (x(t)) \mbox{ }\mbox{for almost every}\mbox{ } t>0
\end{equation*}

\noindent such that
\begin{equation*}\label{strongsol01}
\qquad \dot{x}(t) +  \xi(t) +  \beta (t) \eta (t) =  0.
\end{equation*}

\noindent In particular,  $x(t) \in dom(\partial\Phi) \cap dom(\partial \Psi) \mbox{ }$ for almost every $t>0$.

Existence of strong solutions of nonautonomous monotone differential inclusions is a nontrivial topic. This question is not examined in this paper. We take for granted the existence of such trajectories. The interested reader can consult Br\'ezis \cite{Bre1}, Attouch-Damlamian \cite{AD},
 Kenmochi \cite{Ken} for precise conditions insuring the existence of such solutions. In this paper, we shall be concerned  only  with the asymptotic behavior of the trajectories of the above systems.

\subsection{Key assumption.}
 
We shall prove the convergence property (\ref{MAG3}) under the assumption

\begin{equation*}{(\mathcal H_1)}\label{MAG4}\hfill
\mbox{ }\mbox{ } \mbox{ }\forall p\in R(N_C) \mbox{ } \mbox{ }\mbox{ } \int_{0}^{+\infty} \beta (t) \left[\Psi^* \left(\frac{p}{ \beta (t)}\right) - \sigma_C\left(\frac{p}{ \beta (t)}\right)\right]dt < + \infty.
\end{equation*}

\noindent In $(\mathcal H_1)$ we use classical notions and notations from convex analysis: $\Psi^*$ is the Fenchel conjugate of $\Psi$,

$$ \forall y \in \mathcal H \mbox{ }\mbox{ }\mbox{ }\Psi^* (y) = \sup_{x\in   \mathcal H}  \left\{  \langle  y ,  x \rangle - \Psi (x)\right\},$$

\noindent and $ \sigma_C $ is the support function of $C= \mbox{argmin} \Psi$
 
$$ \forall y \in \mathcal H \mbox{ }\mbox{ }\mbox{ }\sigma_C (y) = \sup_{x\in C}   \langle  y ,  x \rangle.$$

\noindent Note that $ \sigma_C $ is equal to the Fenchel conjugate of ${\delta}_C$, where ${\delta}_C$ is the indicator function of $C$.
$ N_C (x)$ is the (outwards) normal cone to $C$ at $x$.
We denote by $R(N_C)$ the range of $N_C$, i.e., $p\in R(N_C)$ iff
$p\in N_C (x)$ for some $x\in C.$

\vspace{0.2cm}

\textit{Analysis of the condition} $(\mathcal H_1)$:
\vspace{0.1cm}
\begin{itemize}
\item a) Note that $\Psi$ enters in $(MAG)$ only via its subdifferential. Thus it is not a restriction to assume $min_{\mathcal H} \Psi =0$. 
For a function $\Psi$ whose minimum is not equal to zero, one should replace in $(\mathcal H_1)$ and in the corresponding statements $\Psi$ by $\Psi - min_{\mathcal H} \Psi.$

\noindent From $\Psi \leq {\delta}_C$ we get  $\Psi^* \geq \left(\delta_C\right)^* =  \sigma_C$   and $ \Psi^*  - \sigma_C  \geq 0.$
 $(\mathcal H_1)$ means that  the  nonnegative function 

$$ t \mapsto \beta (t) \left[\Psi^* \left(\frac{p}{ \beta (t)}\right) - \sigma_C \left(\frac{p}{ \beta (t)}\right)\right]$$ 

\noindent is integrable on $(0, +\infty)$. It is a growth condition on $\beta(.)$ at infinity which depends only on $\Psi$.

\vspace{0.1cm}

\item b) As an illustration, consider the model situation: $\Psi (z) =  \frac{1}{2} dist^2 (z, C) =  \frac{1}{2} \|.\|^2  \mbox{ }+_e \mbox{ } \delta_C$, where $+_e$ denotes the epigraphical sum (also called inf-convolution). From  general properties of  Fenchel transform 

\begin{center}
$\Psi^* (z) =  \frac{1}{2} \|z\|^2  +  \sigma_C (z)$ and $\Psi^* (z) -   \sigma_C (z)  =  \frac{1}{2} \|z\|^2 $. 
\end{center}

\noindent Hence, in this situation

$$(\mathcal H_1) \Longleftrightarrow  \int_{0}^{+\infty} \frac {1}{\beta (t)} dt < + \infty$$

\noindent which is satisfied for example with $\beta (t) = (1+t)^p , \ p>1$.

\end{itemize}

\subsection{Contents.}

Several results concerning the asymptotic convergence analysis hold true in respect of the more general differential inclusion

\begin{equation}\label{MAG5} 
\dot{x}(t) +  A(x(t)) +  \beta (t) \partial \Psi (x(t))\ni 0
\end{equation}

\noindent with $A$   maximal monotone operator (in particular one may
consider $A=  \partial \Phi$), and without  regularity assumptions on
the function $\beta(.)$ which may present oscillations, discontinuities. 
In  section 2, we prove an ergodic convergence result (theorem 2.1) which holds for (\ref{MAG5}), and which  extends  Baillon-Brezis theorem \cite{BB} to a nonautonomous multiscale setting. 

In section 3, we return to  $(MAG)$ system with $A = \partial \Phi$. By using energetic Liapunov methods, under the additional growth condition on $\beta$, namely 
$\dot{\beta} \leq k\beta,$ 
we prove  an  asymptotic weak convergence result (theorem 3.1). This result can be seen as an extension of  Bruck theorem (\cite{Bruck}).
In section 4, we  revisit the asymptotic analysis of the sytem

\begin{equation}\label{MAG6}
\dot{x}(t) +  \epsilon (t) \partial \Phi (x(t)) +  \partial \Psi 
(x(t))\ni 0 
\end{equation}

\noindent where  $\epsilon (t)$ tends to zero as $t$ goes to $+\infty$ and satisfies $  \int_{0}^{+\infty}\epsilon (t) dt = + \infty  $.  Indeed (\ref{MAG6}) can be  derived from $(MAG)$ 
 by  time rescaling.
In section 5, we show that, in the particular case of inf-compact functions (in particular in the finite dimensional case),  convergence results of section 3 hold without  growth condition $\dot{\beta} \leq k\beta$ .

In  last section 6, applications are given to coupled gradient dynamics. In particular, we consider domain decomposition for elliptic PDE's, and  best response dynamical approach to Nash equilibria for potential games.

\section{With a maximal monotone operator: ergodic convergence results}

In this section, we consider the differential inclusion (\ref{MAG5}) with  a maximal monotone operator $A$.
We call it $(MAMI)$:
  
\begin{equation*}\label{MAG7}
 (MAMI) \qquad \dot{x}(t) +  A(x(t)) +  \beta (t) \partial \Psi (x(t))\ni 0.
\end{equation*} 

 This terminology  stands for
 ``Multiscale Asymptotic Monotone Inclusion'' with a justification similar to (MAG):
  one expects that the nonautonomous multiscale differential inclusion $(MAMI)$ enjoy asymptotic
  properties similar to the autonomous monotone inclusion
 
 \begin{equation*}\label{MAG8}
\dot{x}(t) +  (A +  N_C)(x(t))\ni 0.
\end{equation*}
   
We denote by $S$ the set of equilibria  $S:= (A + N_C)^{-1}(0)$. They are solutions of 
the monotone variational inequality
 
  $$   A(x) +  N_C (x)\ni 0. $$
 
We  first examine $(MAMI)$ with a  general maximal monotone operator $A$, in which case we  prove an ergodic convergence result. Then, in the particular  case of a strongly monotone operator $A$,  we  prove  strong convergence of the trajectories towards the unique equilibrium.

\begin{Theorem} \label{Thergo}
Let

\begin{itemize}

\item $A: \mathcal H \rightarrow 2^{\mathcal H}$  be a general maximal monotone operator.

\item $\Psi: \mathcal H \rightarrow \R^+ \cup \{+ \infty\}$ be a
closed convex proper function, such that $C$ = $\mbox{argmin}\Psi = \Psi ^{-1}(0) \neq \emptyset.$

\end{itemize}

Let us assume that,

\begin{itemize}
\item $(\mathcal H_0)\mbox{ }\mbox{ }$ $A + N_C$ is a maximal monotone operator and $S:= (A + N_C)^{-1}(0)$ is  non empty. 

\item
$(\mathcal H_1)$\mbox{ }\mbox{ }  $ \forall p\in R(N_C),$  $\displaystyle \int_{0}^{+\infty} \beta (t) \left[\Psi^* \left(\frac{p}{ \beta (t)}\right) - \sigma_C \left(\frac{p}{ \beta (t)}\right)\right]dt < + \infty.$

\end{itemize}

 Then, for every strong solution trajectory $x(.)$ of the differential
 inclusion $(MAMI)$:
\begin{eqnarray*}
(i) & \mbox{ (weak ergodic convergence) }&\exists x_{\infty}\in S  \
  \mbox{ such that }  \  w-\lim_{t\rightarrow
  +\infty}\mbox{  } \frac{1}{t} \int_0^t x(s) ds = x_\infty;\\
(ii)&&\forall z\in S,\>\lim_{t\rightarrow +\infty} \|x(t) -z\| \> \mbox{ exists};\\
(iii)&\mbox{ (estimation) }&\int_{0}^{+\infty} \beta (t)\Psi (x(t)) dt < +\infty.
\end{eqnarray*}

\end{Theorem}

By taking $\Psi=0$ in  Theorem~\ref{Thergo}, we recover the
Baillon-Br\'ezis result on the ergodic convergence of semi-groups of contractions in Hilbert spaces (generated by maximal monotone operators).

\begin{corollaire}{\cite[Baillon-Br\'ezis]{BB}}
Let
$A: \mathcal H \rightarrow  2^{\mathcal H}$  be a  maximal monotone operator
 such that $A^{-1}0 \neq \emptyset.$ 
Let $x$ be a strong solution of 

$$ \dot{x}(t) +  A(x(t))\ni 0.$$

Then,
$$w-\lim_{t\rightarrow +\infty}\mbox{  }\mbox{  } \frac{1}{t} \int_0^t x(s) ds = x_\infty  \mbox{ }\mbox{exists with} \mbox{ } x_\infty \in A^{-1}0.$$

\end{corollaire}

\begin{remark}
An elementary example (take $A$ equal to the rotation of angle $\frac{\pi}{2}$ in $\R^2$) shows that  ergodic convergence can happen without convergence. 
Clearly, the same holds true for  $(MAMI)$.
\end{remark}

\subsection{Proof of Theorem~\ref{Thergo}.}
 It relies on an  Opial type argument. Owing to the length of the proof, we decompose it in several lemmas.

\begin{lemma}\label{12092008_1} For every $z\in S= (A + N_C)^{-1}(0), \mbox{ } lim_{t \to  +\infty} \parallel x(t)- z \parallel$  exists.
\end{lemma}

{\bf Proof of Lemma~\ref{12092008_1}.} Set 
$$
h_z (t) = \frac{1}{2}\parallel x(t) - z \parallel ^2
$$ 
and compute
the time derivative of $h_z (t)$.  For almost every $t>0$
\begin{eqnarray*}
\dot{h}_z (t) & = & \langle  x(t) -z,  \dot{x}(t) \rangle\\
 & = & \langle  x(t) -z,  - \xi(t)  -  \beta (t) \eta (t) \rangle,
\end{eqnarray*}

\noindent where, by definition of $x(.)$ solution of $(MAMI)$,

$$
\qquad \dot{x}(t) +  \xi(t) +  \beta (t) \eta (t) =  0
$$
with
$$
\xi(t) \in A(x(t)) \mbox{ }\mbox{and}\mbox{ } \eta (t) \in \partial
\Psi (x(t)).
$$
Equivalently,

\begin{equation}\label{equ00}
\dot{h}_z (t) +  \langle  \xi(t),  x(t) -z \rangle  + \beta (t) \langle \eta (t), x(t) -z \rangle  =  0 .
\end{equation}

\noindent Since $z\in S \subset C$ and $\eta (t) \in \partial \Psi (x(t)),$  we have

$$
 0 = \Psi (z) \geq  \Psi(x(t))  +  \langle \eta (t), z - x(t) \rangle,  
$$

that is, 

\begin{equation}\label{equ21}
 \langle \eta (t), x(t) -  z \rangle  \geq   \Psi(x(t)). 
\end{equation}

Since $z\in S$,  we have $Az + N_C (z) \ni 0$  \
i.e., there exists some $p\in N_C (z)$ such that $-p \in Az.$ 
By monotonicity of $A$ and $\xi(t) \in A(x(t))$ 
 
\begin{equation}\label{equ22}
\langle  \xi(t),  x(t) -z \rangle  \geq \langle  -p,  x(t) -z \rangle. 
\end{equation}

In view of (\ref{equ00}), (\ref{equ21}) and (\ref{equ22}), and since $\beta(t) > 0$ we obtain

\begin{equation*}\label{eq_base}
\dot{h}_z (t) \leq   \langle  p,  x(t) -z \rangle  - \beta (t) \Psi
(x(t)).
\end{equation*}

Let us rewrite this inequality as

\begin{equation}\label{equ3}
\dot{h}_z (t) \leq  \beta (t) \left[ \left\langle  \frac{p}{ \beta (t)},  x(t) -z \right\rangle  -  \Psi (x(t))\right].
\end{equation}

Since we have no prior information on $x(t)$, let us take the supremum of this last expression with respect to $x$

$$
\dot{h}_z (t) \leq   \beta (t) \left[ \sup_{x\in\mathcal H} \left\{ \left\langle  \frac{p}{ \beta (t)},  x \right\rangle  -  \Psi (x)\right\} 
-  \left\langle \frac{p}{\beta(t)}, z\right\rangle \right],
$$

which makes appear the Fenchel conjugate of $\Psi$, namely ${\Psi}^*$, and gives 

\begin{equation}\label{equ4}
\dot{h}_z (t) \leq   \beta (t) \left[\Psi^*\left(\frac{p}{ \beta (t)}\right) - \left\langle \frac{p}{ \beta (t)}, z \right\rangle  \right].
\end{equation}

Let us now examine the term $\langle \frac{p}{ \beta (t)}, z \rangle $ in (\ref{equ4}).
 Since  $z\in C$ and $p\in N_C (z)$ we have

$$
 \forall x \in C   \mbox{  }\mbox{  } \left\langle p, x - z \right\rangle  \leq 0
$$

which implies

\begin{eqnarray*}
\langle p, z \rangle  & = & \sup_{x\in C}  \langle p, x \rangle 
 =  {\sigma}_C (p).
\end{eqnarray*}

By positive homogeneity of ${\sigma}_C$ and $\beta (t)>0$, we obtain

\begin{equation}\label{equ5}
\left\langle \frac{p}{ \beta (t)}, z \right\rangle = {\sigma}_C \left(\frac{p}{ \beta (t)}\right).
\end{equation}

Collecting (\ref {equ4}, \ref{equ5}) we finally obtain

\begin{equation}\label{equ6}
\dot{h}_z (t) \leq   \beta (t) \left[\Psi^*  \left(\frac{p}{ \beta (t)}\right) - {\sigma}_C  \left(\frac{p}{ \beta (t)}\right)  \right].
\end{equation}

Note that $\Psi \leq {\delta}_C$, hence ${\Psi}^* \geq {\sigma}_C$ and $ {\Psi}^*  - {\sigma}_C  \geq 0.$

It follows from  (\ref{equ6})  and assumption $(\mathcal H_1)$ that ${(\dot{h}_z)}_+  \in L^1 (0, + \infty)$,
which classically implies that $\mbox{lim}_{t\to +\infty} h_z (t)$
exists in $\R$.$\Box$

\vspace{0.3cm}

\begin{lemma}\label{12092008_2}  For each $t > 0$ set
$$ 
X(t)  =   \frac{1}{t}  \int_0^t  x(s)ds . 
$$
Every weak limit point of $X(.)$ belongs to $S$. 
\end{lemma}

{\bf Proof of Lemma~\ref{12092008_2}.}
Let $t_n \to +\infty$ and suppose $X(t_n) \rightharpoonup X_{\infty}$
(weak convergence in $\mathcal H$).
Take an arbitrary $z\in C \cap \mbox{dom} A$ , and $y\in  (A + N_C)(z)$.
Consider again the
function $h_z (t) = \frac{1}{2}\parallel x(t) - z \parallel ^2$.
Recall  (\ref{equ00}):
$$
\dot{h}_z (t) +  \langle  \xi(t),  x(t) -z \rangle  + \beta (t)
\langle \eta (t), x(t) -z \rangle  =  0 .
$$

For some $p\in N_C (z)$, $y-p\in A(z)$. By monotonicity of $A$ and $ \xi(t) \in A(x(t))$

$$
\langle  \xi(t),  x(t) -z \rangle  \geq \langle  y-p,  x(t) -z \rangle.
$$

Recall (\ref{equ21})

$$   
\langle \eta (t), x(t) -  z \rangle  \geq   \Psi(x(t)) 
$$

and we obtain 

$$
\dot{h}_z (t) +  \langle  y ,  x(t) -z \rangle  \leq  \beta (t)
\left[ \left\langle \frac{p}{ \beta (t)}, x(t) -z \right\rangle   - \Psi(x(t))
\right].
$$

Hence

\begin{equation*}\label{equ9}
\dot{h}_z (t) +  \langle  y ,  x(t) -z \rangle  \leq  \beta (t) \left[\Psi^* \left(\frac{p}{ \beta (t)}\right) - {\sigma}_C \left(\frac{p}{ \beta (t)}\right)  \right].                    
\end{equation*}

Let us integrate from $0$ to $t$

\begin{equation*}
 h_z (t) +  \left\langle  y , \int_0^t  x(s)ds -tz \right\rangle  \leq h_z (0) +  \int_0^t \beta (s) \left[\Psi^* \left(\frac{p}{ \beta (s)}\right) - {\sigma}_C \left(\frac{p}{ \beta (s)}\right)  \right] ds.
\end{equation*}

After division by $t$, and taking account of $h_z \geq 0$, one obtains

\begin{eqnarray*}
\langle  y , X(t) - z \rangle  & \leq  & \frac{1}{t}  h_z (0) +   \frac{1}{t} \int_0^t \beta (s) \left[\Psi^* \left(\frac{p}{ \beta (s)}\right) - {\sigma}_C \left(\frac{p}{ \beta (s)}\right)  \right] ds\\
 & \leq & \frac{1}{t}  h_z (0) +   \frac{1}{t} \int_0^{+\infty} \beta (s) \left[\Psi^* \left(\frac{p}{ \beta (s)}\right) - {\sigma}_C \left(\frac{p}{ \beta (s)}\right)  \right] ds\\
  & \leq & \frac{c}{t}
\end{eqnarray*}

where

$$
 c := h_z (0) + \int_0^{+\infty} \beta (s) \left[\Psi^* \left(\frac{p}{ \beta (s)}\right) - {\sigma}_C \left(\frac{p}{ \beta (s)}\right)  \right] ds,
$$ 

is,  by assumption $(\mathcal H_1)$,  a finite positive number.

Recall that  $t_n \to +\infty$ and  $X(t_n) \rightharpoonup X_{\infty}$  (weak convergence in $\mathcal H$).
Passing to the limit as $ n\to +\infty$ on 

$$\langle  y , X(t_n) - z \rangle \leq  \frac{c}{t_n}$$

we finally obtain

$$\langle  y , X_{\infty} - z \rangle \leq  0.$$

An equivalent formulation is

\begin{equation*}\label{equ10}
   \langle 0 - y, X_{\infty} - z \rangle   \geq 0.                         
\end{equation*}

The inequality  being true for any $y\in \mbox{dom}( A + N_C)$ and any $y\in (A +
N_C)z$,  by maximal monotonicity of the  operator $A + N_C$  (Assumption
$(\mathcal H_0)$),  we obtain $0 \in (A + N_C)(X_{\infty})$, that is 
$X_{\infty} \in S.\Box$

\vspace{0.3cm}

Just like in Passty \cite{Passty}, we conclude to the ergodic
convergence of the trajectories thanks to  Lemmas~\ref{12092008_1} and~\ref{12092008_2} and the following ergodic variant of Opial's lemma~\cite{opial}.

\begin{lemma}\label{passty} Let  $\mathcal H$ be a Hilbert space, $S$
a non empty subset of $\mathcal H$ and  $x: [0, +\infty) \to \mathcal
H$ a map. Set
$$
X(t)  =   \frac{1}{t}  \int_0^t  x(s)ds
$$
and assume that
\begin{eqnarray*}
(i) &&\mbox{for every }z\in S,\>  \lim_{t \to  +\infty} \parallel
x(t)- z \parallel  \mbox{ exists};\\
(ii) &&\mbox{every weak limit point of the map }X \mbox{ belongs to }S. 
\end{eqnarray*}
Then
 $$
w-\lim_{t \to  +\infty} X(t) = X_{\infty} \mbox{ for some element
}X_{\infty}\in S. 
$$
\end{lemma}

{\bf Proof of Lemma~\ref{passty}.} First note that $x$ is bounded
(from $(i)$), thus $X$ also and it is sufficient to prove  uniqueness
of weak limit points.
Let  $X(t_{n_1}) \rightharpoonup X_{\infty, 1}$ and $X(t_{n_2}) \rightharpoonup X_{\infty, 1}$
be two weak converging subsequences. \\
Let us prove  that $X_{\infty, 1} = X_{\infty, 2}$.
By $(ii)$, $X_{\infty, 1}$ and  $ X_{\infty, 2}$ belong to $S$. Hence, 
by $(i)$,  $\mbox{   lim}_{t \to  +\infty} \parallel x(t)- X_{\infty, 1} \parallel^2 $  and $\mbox{   lim}_{t \to  +\infty} \parallel x(t)- X_{\infty, 2} \parallel^2 $ exist. As a consequence, the following limit exist

\begin{center}
 $\mbox{   lim}_{t \to  +\infty} \left[\parallel x(t)- X_{\infty, 1} \parallel^2 - \parallel x(t)- X_{\infty, 2} \parallel^2 \right].$
\end{center}

After simplification

$$
\lim_{t \to  +\infty}  \langle  x(t) ,  X_{\infty, 2} -  X_{\infty,
1} \rangle  \ \mbox{ exists.}
$$

As a general classical result (Cesaro),  convergence implies  ergodic convergence. Hence

$$
\lim_{t \to  +\infty}  \langle  X(t) ,  X_{\infty, 2} -  X_{\infty, 1} \rangle  \  \mbox{ exists.} 
$$

In particular, 

$$
 \lim_{t_{n_1} \to  +\infty}  \langle  X(t_{n_1}) ,  X_{\infty, 2} -  X_{\infty, 1} \rangle = 
 \lim_{t_{n_2} \to  +\infty}  \langle  X(t_{n_2}) ,  X_{\infty, 2} -  X_{\infty, 1} \rangle,
$$

that is

\begin{center}
 $  \langle  X_{\infty, 1} ,  X_{\infty, 2} -  X_{\infty, 1} \rangle = \langle X_{\infty, 2} ,  X_{\infty, 2} -  X_{\infty, 1} \rangle$  
\end{center}

and   $\parallel X_{\infty, 2} -  X_{\infty, 1} \parallel^2  =  0$, \ which ends the proof of the lemma.$\Box$

\vspace{0.3cm}

\textit{End of the proof of theorem \ref{Thergo}}

Let us complete the proof of  theorem  \ref{Thergo} and prove  the estimation

$$
\int_{0}^{+\infty} \beta (t)\Psi (x(t)) dt < +\infty.
$$

Let us return to (\ref{equ3})

$$
\dot{h}_z (t) \leq  \beta (t) \left[ \left\langle  \frac{p}{ \beta (t)},  x(t) -z \right\rangle  -  \Psi (x(t))\right]
$$

which can be written in a splitted form as

$$
\dot{h}_z (t) +  \frac{\beta (t)}{2}\Psi (x(t))  \leq  \frac{\beta (t)}{2} \left[\left \langle  \frac{2p}{ \beta (t)},  x(t) -z \right\rangle  -  \Psi (x(t))\right].
$$

By using a  device similar to the proof of lemma~\ref{12092008_1}

\begin{eqnarray*}
 \dot{h}_z (t) +  \frac{\beta (t)}{2}\Psi (x(t)) & \leq  & \frac{\beta (t)}{2} \left[ \sup_{x\in\mathcal H} \left\{ \left\langle  \frac{2p}{ \beta (t)},  x \right\rangle  -  \Psi (x)\right\} 
-  \left\langle \frac{2p}{\beta(t)}, z\right\rangle \right]\\
 & \leq & \frac{\beta (t)}{2} \left[ {\Psi}^* \left( \frac{2p}{ \beta (t)}\right) -  {\sigma}_C \left(\frac{2p}{ \beta (t)}\right)   \right].\\
\end{eqnarray*}

Let us integrate this last inequality from $0$ to $\tau$

\begin{equation}\label{equ11}
h_z (\tau) - h_z (0) +  \frac{1}{2}\int_0^{\tau} \beta (t) \Psi (x(t)) dt \leq \frac{1}{2} \int_{0}^{+\infty} \beta (t) \left[\Psi^* \left(\frac{2p}{ \beta (t)}\right) - \sigma_C \left(\frac{2p}{ \beta (t)}\right)\right]dt.
\end{equation}

By assumption $(\mathcal H_1)$, the second member of (\ref{equ11}) is a finite quantity (note that $p\in R( N_C) $ implies $2p\in R(N_C) $). 

This being true for any $\tau > 0$, we finally obtain 
$$
 \int_{0}^{+\infty} \beta (t)\Psi (x(t)) dt < +\infty.
$$

\subsection{The case $A$ strongly monotone}

The argument developed in the preceding section allows to conclude  to  convergence of the trajectories when the operator $A: \mathcal H \rightarrow \mathcal H$  satisfies a strong monotonicity property. We recall that $A$ is said to be strongly monotone if there exists some $\alpha >0$ such that for any $x\in domA$, $y\in domA$ and any $\xi \in Ax$, $\eta \in Ay$,

$$ \langle  \xi - \eta,  x -y \rangle  \geq \alpha \parallel x-y\parallel^2. $$

\begin{Theorem} \label{Thsm} Let us assume that,

\begin{itemize}
\item $(\mathcal H_0)\mbox{ }\mbox{ }$ $A + N_C$ is a maximal monotone operator and $S:= (A + N_C)^{-1}(0)$ is  non empty. 

\item
$(\mathcal H_1)$\mbox{ }\mbox{ }  $ \forall p\in R(N_C)  \mbox{ } \mbox{ }\mbox{ } \displaystyle \int_{0}^{+\infty} \beta (t) \left[\Psi^* \left(\frac{p}{ \beta (t)}\right) - \sigma_C \left(\frac{p}{ \beta (t)}\right)\right]dt < + \infty.$

\item
$(\mathcal H_2)\mbox{ }\mbox{ }$  $A$ is a strongly monotone operator.

\end{itemize}

 Then, there exists a unique equilibrium $\xb$, i.,e., $S = \{\xb\}$ with $A\xb + N_C (\xb) \ni 0$ and
any strong solution trajectory $x(.)$ of the differential inclusion $(MAMI)$ strongly converges to $\xb$ as 
$t\rightarrow +\infty$:
 
$$ s-lim_{t\rightarrow +\infty}\mbox{  }\mbox{  } x(t)  = \xb. $$

\end{Theorem}

\textbf{{Proof of Theorem}~\ref{Thsm}}. By strong monotonicity of $A$, the operator $\mbox{  } A + N_C\mbox{  }$ is also strongly monotone.
As a consequence, there exists a unique solution $\xb$ to the inclusion

\begin{center}
$A\xb + N_C (\xb) \ni 0$
\end{center}

and there exists some $p \in N_C(\xb)$ such that $A\xb \ni -p.$

Let us return to (\ref{equ00}), with $z = \xb$ and $h_{\xb}(t) = \frac{1}{2}\parallel x(t) -\xb \parallel ^2$,

\begin{equation}\label{equ01}
\dot{h}_{\xb} (t) + \langle  \xi (t),  x(t) - \xb \rangle  + \beta (t) \langle \eta (t), x(t) - \xb \rangle  =  0 .
\end{equation}

Let us rewrite (\ref{equ01}) as 

$$
\dot{h}_{\xb} (t) + \langle  \xi (t) - (-p),  x(t) - \xb \rangle  - \langle p,  x(t) - \xb \rangle + \beta (t) \langle  \eta (t), x(t) - \xb \rangle  =  0,
$$

and use the strong monotonicity of $A$ together with $ \mbox{  } \xi
(t) \in A (x(t))$ to obtain, for some $\alpha>0$,

$$
\dot{h}_{\xb} (t) + \alpha \parallel x(t) -\xb \parallel^2   - \langle p,  x(t) - \xb \rangle + \beta (t) \langle \eta (t), x(t) - \xb \rangle  \leq  0.
$$

By using the convex subdifferential inequality

$$  
0 = \Psi (\xb) \geq  \Psi(x(t))  +  \langle \eta (t), \xb - x(t) \rangle  
$$

we deduce

$$
\dot{h}_{\xb} (t) + \alpha \parallel x(t) -\xb \parallel^2   - \langle p,  x(t) - \xb \rangle + \beta (t) \Psi(x(t))  \leq  0.
$$

Hence,

$$
\dot{h}_{\xb} (t) + \alpha \parallel x(t) -\xb \parallel^2  \leq   \beta (t) \left[ \left\langle \frac{p}{\beta (t)},  x(t) - \xb \right\rangle -    \Psi(x(t))\right].
$$

By using a  device similar to theorem \ref{Thergo}

$$
\dot{h}_{\xb} (t) + \alpha \parallel x(t) -\xb \parallel^2  \leq   \beta (t) \left[ {\Psi}^* \left(\frac{p}{\beta (t)}\right)    
- {\sigma}_C \left(\frac{p}{\beta (t)}\right) \right].
$$

After integration of this inequality from $0$ to $t$,

$$h_{\xb} (t) + \alpha \int_0^t  \parallel x(s) -\xb \parallel^2  ds  \leq  h_{\xb} (0) +  \int_0^{+\infty} \beta (s)\left[ {\Psi}^* \left(\frac{p}{\beta (s)}\right)    
- {\sigma}_C \left(\frac{p}{\beta (s)}\right) \right]ds,
$$

 by using assumption $(\mathcal H_1)$, one obtains

$$
 \int_0^{+\infty} \parallel x(t) -\xb \parallel^2  dt < +\infty.
$$

On the other hand, by  theorem \ref{Thergo},  $\mbox{lim}_{t\to +\infty} \parallel x(t) -\xb \parallel$
exists. As a consequence, this limit is  equal to zero, i.e. $\mbox{lim}_{t\to +\infty} \parallel x(t) -\xb \parallel = 0$, which ends the proof of theorem \ref{Thsm}. $\Box$

\section{The subdifferential case: weak convergence results}\label{section_sub}

In this section, we consider the dynamical system

$$(MAG) \qquad \dot{x}(t) +  \partial \Phi (x(t)) +  \beta (t)
\partial \Psi (x(t))\ni 0.$$

\noindent When $\Psi=0$, $(MAG)$ boils down to the classical steepest descent differential inclusion
$$\dot{x}(t) +  \partial \Phi (x(t))\ni 0,$$
studied by Br\'ezis ~\cite{Bre1}~\cite{Bre2}, Bruck~\cite{Bruck} and Baillon~\cite{Baillon}. In accordance with these studies,  in our main
result, we are going to show that each trajectory of $(MAG)$  weakly convergences  to a
minimizer of $\Psi$, which also minimizes $\Phi$ over all minima of
$\Psi$.

Before stating  our  result precisely, let us specify the notion of solution. We recall that a
solution $x(.)$ of $(MAG)$ is continuous on $\left[0, + \infty \right)$, absolutely continuous on any bounded interval
 $\left[0,T\right]$ with $T< +\infty$,
 and it it assumed that the equation holds almost
 everywhere.  Moreover,  in this section, we will assume that there exist two functions $\xi(.)$ and $\eta (.)$ which are integrable on any 
bounded interval $\left[0,T\right]$ and such that

\begin{equation}\label{strongsol1}
\qquad \dot{x}(t) +  \xi(t) +  \beta (t) \eta (t) =  0
\end{equation}

with 

\begin{equation}\label{strongsol2}
\xi(t) \in  \partial \Phi(x(t)) \mbox{ }\mbox{and}\mbox{ } \eta (t) \in \partial \Psi (x(t)) \mbox{ }\mbox{for almost every}\mbox{ } t>0.
\end{equation}

In several cases, for example if $\Psi$ is differentiable, with a
Lipschitz gradient $\nabla \Psi$ on the bounded sets, or if the
subdifferentials $\partial \Phi$ and  $\partial \Psi$ satisfy an angle
condition, the assumption that the two functions $\xi(.)$ and $\eta
(.)$ are locally integrable is automatically fulfilled, provided that
the trajectory $x$ is  absolutely continuous. However, in general, it is a
nontrivial issue for which we refer to Attouch-Damlamian \cite{AD}.
We will denote by $S=\mbox{argmin}\{ \Phi|\mbox{argmin}\Psi\}$ the set of equilibria, which is a dense notation for 
$S = \{z\in \mbox{argmin}\Psi: \Phi(z) \leq \Phi (x) \  \mbox{ for all }  x\in \mbox{argmin}\Psi\}$.

\begin{Theorem} \label{Theweak}
Let

\begin{itemize}
\item $\Psi: \mathcal H \rightarrow \R^+ \cup \{+ \infty\}$ be a
closed convex proper function, such that $C$ = $\mbox{argmin}\Psi = \Psi ^{-1}(0) \neq \emptyset.$
\item $\Phi: \mathcal H \rightarrow \R^{\phantom{+}} \cup \{+ \infty\}$ be a
closed convex proper function, such that
$S=\mbox{argmin}\{ \Phi|\mbox{argmin}\Psi\} \neq \emptyset.$
\end{itemize}

Let us assume that,

\begin{itemize}
\item
$(\mathcal H_1)$\mbox{ }\mbox{ }  $ \forall p\in R(N_C),$
  $\displaystyle \int_{0}^{+\infty} \beta (t) \left[\Psi^*
  \left(\frac{p}{ \beta (t)}\right) - \sigma_C \left(\frac{p}{ \beta
    (t)}\right)\right]dt < + \infty.$

\item
$(\mathcal H_2)$\mbox{ }\mbox{ } $\beta: \R^+ \rightarrow \R^+$  is a function  of class $C^1$,
  such that $\lim_{t\to+\infty} \beta (t)=+\infty$, and for some $k\geq 0$ and $t_0 \geq 0$
\begin{center}
 $0\leq\dot{\beta}(t)\leq k \beta(t)$ for all $t\geq t_0$.
\end{center}

\end{itemize}

Let $x$ be a strong solution of $(MAG)$. Then:

\begin{eqnarray*}
(i)&\mbox{ weak convergence }&\exists x_{\infty}\in S=\mbox{argmin}\{
  \Phi|\mbox{argmin}\Psi\},\qquad  w-\lim_{t\to+\infty}
  x(t)=x_{\infty};\\
(ii)&\mbox{ minimizing properties } &\lim_{t\to+\infty} \Psi(x(t))=
  0;\\
&&\lim_{t\to +\infty} \Phi(x(t)) = \min\Phi |_{\mbox{argmin}\Psi};\\
(iii)&&\forall z\in S  \  \lim_{t\to+\infty} \|x(t) -z\| \mbox{  exists  };\\
(iv)&\mbox{ estimations }&\lim_{t\to +\infty} \beta(t)\Psi(x(t))=0;\\
&&  \int_{0}^{+\infty} \beta (t)\Psi (x(t)) dt < +\infty;\\
&&\limsup_{\tau\to +\infty} \int_{0}^{\tau}\Phi (x(t))- \min\Phi |_{\mbox{argmin}\Psi}  dt < +\infty.
\end{eqnarray*}

\end{Theorem}

By taking $\Psi=0$ in  Theorem~\ref{Theweak}, we recover the
convergence result of Bruck on the steepest descent method.

\begin{corollaire}{\cite[Bruck, Theorem 4]{Bruck}}
Let
$\Phi: \mathcal H \rightarrow \R^+ \cup \{+ \infty\}$ be a
closed convex proper function, such that $\mbox{argmin} \Phi \neq
\emptyset.$ 
Let $x$ be a strong solution of 
$$(SD)\qquad \dot{x}(t) +  \partial \Phi (x(t))\ni 0.$$
Then $x$ weakly converges to a point in $\mbox{argmin}
\Phi$.
\end{corollaire}

\begin{remark}
The counterexample of Baillon~\cite{Baillon} shows that one may not
have strong convergence for $(SD)$. Of course, the same holds for
$(MAG)$.
\end{remark}

\textbf{Proof of Theorem~\ref{Theweak}.}
As in Bruck~\cite{Bruck}, the weak convergence is consequence of
Opial's lemma, after showing the convergence of $ \|x(.) -z\|$ for
every $z\in S$, and that every weak limit point of $x$ belongs to
$S$. The proof is not short, and we decompose it in several lemmas.
 For an element $z\in S=\mbox{argmin}\{
  \Phi|\mbox{argmin}\Psi\}$, we define the function $h_z:\R_+\to \R_+$
  by
$$
h_z(t)=\frac{1}{2}\|x(t) - z \| ^2.
$$

We give the estimations on the function $h_z$, that we will use
in the proof. Compute the derivative of $h_z$ and use the system
equation $(MAG)$, $(\ref{strongsol1})$ and $(\ref{strongsol2})$, for
a.e. $t$:
\begin{eqnarray*}
\dot{h}_z(t)&=&\langle \dot{x}(t), x(t)-z\rangle\\
&=&\langle -\xi(t)-\beta(t)\eta(t), x(t)-z\rangle,
\end{eqnarray*}
with $\xi(t)\in \partial  \Phi(x(t))$ and  $\eta (t) \in \partial \Psi
(x(t))$. Thus, using $z\in\mbox{argmin}\Psi$,
\begin{eqnarray*}
\Phi(z)&\geq& \Phi(x(t))+\langle \xi(t), z - x(t)\rangle;\\
0=\Psi(z)&\geq& \Psi(x(t))+\langle \eta(t), z - x(t)\rangle.
\end{eqnarray*}
We deduce (note that $\beta(t) > 0$)
\begin{equation*}\label{eq_base_0}
\dot{h}_z(t)+ \Phi(x(t))- \Phi(z)+\beta(t)\Psi(x(t))\leq 0.
\end{equation*}
Since $z\in S=\mbox{argmin}\{ \Phi|\mbox{argmin}\Psi\}$, write the
  first order necessary condition (with $C = \mbox{argmin}\Psi$)
$$
0\in \partial \Phi(z)+N_C(z),
$$
and there exists $p\in N_C(z)$ such that $-p \in \partial \Phi(z)$.
Thus 
\begin{eqnarray*}
\Phi(x(t))&\geq& \Phi(z   )+\langle -p, x(t)-z\rangle
\end{eqnarray*}
and
\begin{equation}\label{eq_base_1.5}
\dot{h}_z (t) + \beta (t) \Psi
(x(t))+\langle  -p,  x(t) -z \rangle\leq \dot{h}_z(t)+ \Phi(x(t))- \Phi(z)+\beta(t)\Psi(x(t))\leq 0.
\end{equation}

\begin{lemma}\label{19_1_2009_1}  For every $z\in S$,  as $t\to +\infty$,
\begin{itemize}
\item
(i) $ \|x(t)
  -z\|$ converges in $\R$;
\item (ii) $t\mapsto
\int_0^t \Phi(x(s))-\Phi(z) +\beta (s) \Psi
(x(s)) ds$  converges in $\R$;
\item
(iii)  $\int_0^t \langle p, x(s)-z\rangle ds$
    converges in $\R$;
\item
(iv)  Moreover, $\int_0^{+\infty} \beta(t) \Psi (x(t)) dt< +\infty.$
\end{itemize}

\end{lemma}

{\bf Proof of Lemma~\ref{19_1_2009_1}.} Parts (i) and (iv) were
already proved in Section 2. However, they  are
obtained easily with the others, which permits an easier reading.
Recall that, by definition of the Fenchel conjugate ${\Psi}^*$ of
$\Psi$, 
$$
\Psi^*\left(\frac{p}{ \beta (t)}\right) \geq \left\langle  \frac{p}{
  \beta (t)},  x(t) \right\rangle  -  \Psi (x(t)),
$$
and that $z\in C$ and $p\in N_C (z)$ imply
$$
{\sigma}_C \left(\frac{p}{ \beta (t)}\right)=\left\langle  \frac{p}{
  \beta (t)},  z \right\rangle,
$$
thus, in view of  (\ref{eq_base_1.5})
$$
\dot{h}_z (t) + \beta (t)\left(-\Psi^*  \left(\frac{p}{ \beta
  (t)}\right) + {\sigma}_C  \left(\frac{p}{ \beta (t)}\right)
\right)\leq 
\dot{h}_z (t) + \beta (t) \Psi
(x(t))+\langle  -p,  x(t) -z \rangle \leq 0.
$$
From Assumption $(\mathcal H_1)$, by integrating the above equation
between two large enough real numbers, we first deduce that the function
$h_z$ satisfies the Cauchy criterion, and is bounded, thus converges in
$\R$, that's $(i)$. Then we deduce the convergence of the function
$$
t\mapsto \int_0^t \beta (s) \Psi
(x(s))+\langle  -p,  x(s) -z \rangle ds,
$$
and, in view of (\ref{eq_base_1.5}),
 of the function
$$
t\mapsto \int_0^t \Phi(x(s))-\Phi(z) +\beta (s) \Psi
(x(s)) ds,
$$
that's $(ii)$.
Since $\Psi \geq 0$, and  in view of (\ref{eq_base_1.5}),
\begin{equation*}\label{19_1_2009_2}
\dot{h}_z (t) + \frac{\beta (t)}{2} \Psi
(x(t))+\langle  -p,  x(t) -z \rangle \leq 
\dot{h}_z (t) + \beta (t) \Psi
(x(t))+\langle  -p,  x(t) -z \rangle \leq 0.
\end{equation*}
As 
$\frac{2p}{\beta(t)}\in N_C (z)$,
$$
\dot{h}_z (t) +  \frac{\beta (t)}{2}\left(-\Psi^*  \left(\frac{ 2 p}{ \beta
  (t)}\right) + {\sigma}_C  \left(\frac{2 p}{ \beta (t)}\right)
\right)\leq 
\dot{h}_z (t) +  \frac{\beta (t)}{2} \Psi
(x(t))+\langle  -p,  x(t) -z \rangle \leq 0,
$$
which, in view of $(\mathcal H_1)$, implies the convergence of 
$$
t\mapsto \int_0^t \frac{\beta (s)}{2} \Psi
(x(s))+\langle  -p,  x(s) -z \rangle ds.
$$
From $(ii)$   we deduce that
$$
\int_0^{+\infty} \beta (t) \Psi (x(t))dt < +\infty,
$$
and that $t\mapsto\int_0^t \langle p, x(s)-z\rangle ds$
    converges in $\R$, these are respectively $(iv)$ and $(iii)$.

Let us recall and translate a lemma from Br\'ezis~\cite{Bre1}, in the
special case of the subdifferential of a convex function:
 
\begin{lemma}{\cite[Lemme 4, p73]{Bre1}} \label{lemma_Brezis} Let  $\Phi: \mathcal H \rightarrow \R \cup \{+\infty\}$ be a closed convex proper function. Let
  $x\in L^2(0, T; H)$ be such that $\dot{x}\in L^2 (0, T; H)$ and
  $x(t)\in Dom(\partial \Phi)$ a.e. $t$. Assume that there exists
  $\xi\in L^2(0, T; H)$ such that $\xi(t)\in \partial \Phi(x(t))$ for
  a.e. $t$. Then the function $t\mapsto \Phi(x(t))$ is absolutely
  continuous and for every $t$ such that $x(t)\in Dom(\partial \Phi)$,
  $x$ and  $\Phi(x)$ are differentiable at $t$, we have
$$
\forall h\in  \partial \Phi(x(t)),\qquad \frac{d}{dt}  \Phi(x(t))=\langle \dot{x}(t), h\rangle.
$$

\end{lemma}
\begin{lemma} \label{20_1_2009_1}
Let 
$$
E_1(t):= \frac{\Phi(x(t))}{\beta(t)}+\Psi(x(t)).
$$
Then the function $E_1$ is  absolutely continuous, and, for a.e. $t$,
$$
\dot{E}_1(t)= -\frac{\dot{\beta}(t)}{\beta(t)^2}\Phi(x(t))
-\frac{|\dot{x}(t)|^2}{\beta(t)}.
$$
The function $E_1$ converges to zero, and $\lim_{t\to +\infty}\Psi(x(t))=0$.

\end{lemma}

{\bf Proof of Lemma~\ref{20_1_2009_1}.} 
From Lemma~\ref{lemma_Brezis},
the functions  $t\mapsto \Phi(x(t))$ and $t\mapsto \Psi(x(t))$ are
absolutely continuous, and for a.e. $t$,
\begin{eqnarray*}
\frac{d}{dt}  \Phi(x(t))&=&\langle \dot{x}(t), \xi(t)\rangle\\
\frac{d}{dt}  \Psi(x(t))&=&\langle \dot{x}(t), \eta(t)\rangle.
\end{eqnarray*}
Thus the function $E_1$ is absolutely continuous (recall that $\beta$ is of class $C^1$) and,  by using (\ref{strongsol1}), for a.e. t,
\begin{eqnarray*}
\dot{E}_1(t) &=& \frac{1}{\beta(t)}\langle \dot{x}(t),
\xi(t)\rangle-\frac{\dot{\beta}(t)}{\beta(t)^2}\Phi(x(t))+\langle
\dot{x}(t), \eta(t)\rangle\\
&=&-\frac{\dot{\beta}(t)}{\beta(t)^2}\Phi(x(t)) -\frac{|\dot{x}(t)|^2}{\beta(t)}.
\end{eqnarray*}

\noindent By Lemma~\ref{19_1_2009_1}, the trajectory
$x$ is bounded. Write, for example,
$$
\Phi(x(t))\geq \Phi(z   )+\langle -p, x-z\rangle,
$$
to see that $\Phi(x(t))$ is bounded from below.
By Assumption $(\mathcal H_2)$, the function $\beta$ is non
decreasing. Hence
\begin{eqnarray*}
\dot{E}_1(t) &\leq& -\frac{\dot{\beta}(t)}{\beta(t)^2} \inf_{s\geq 0} \Phi(x(s)) -\frac{|\dot{x}(t)|^2}{\beta(t)}
\end{eqnarray*}
which implies the convergence of  the function $E_1$, recalling that
$\lim_{t\to +\infty}\beta(t)=+\infty$.
From Lemma~\ref{19_1_2009_1}, $(ii)$, we have 
$$
\liminf_{t\to +\infty}\Phi(x(t))-\Phi(z) +\beta (t) \Psi
(x(t))\leq 0.
$$
Using that $\Phi(x(t))$ is bounded from
below and that $ \Psi (x(t))\geq 0$, we see that the above lower limit belongs to $\R$. As a consequence, we can take a sequence $(t_n)$ which tends to $+\infty$
such that $\Phi(x(t_n))-\Phi(z) +\beta (t_n) \Psi
(x(t_n))$ converges, and obtain 
$$
\lim_{t\to +\infty} E_1(t)=\lim_{n\to +\infty} E_1(t_n)=\lim_{n\to +\infty} \frac{1}{\beta(t_n)}\left(\Phi(x(t_n))-\Phi(z) +\beta (t_n) \Psi
(x(t_n))\right)+ \frac{1}{\beta(t_n)}\Phi(z)=0.
$$
\noindent Write
$$
0\leq  \Psi(x(t))\leq  \Psi(x(t))+ \frac{1}{\beta(t)}\left(\Phi(x(t))-
\inf_{s\geq 0} \Phi(x(s))\right) = E_1 (t) - \frac{1}{\beta(t)}\inf_{s\geq 0} \Phi(x(s))
$$

\noindent to deduce
$$
\lim_{t\to +\infty} \Psi(x(t))=0.
$$

\begin{lemma}\label{20_1_2009_2} 
\begin{eqnarray*}
& &\liminf_{t\to +\infty} \langle  -p,  x(t) -z \rangle=0;\\
& &\liminf_{t\to +\infty}\Phi(x(t))\geq \Phi(z).
\end{eqnarray*}
\end{lemma}

{\bf Proof of Lemma~\ref{20_1_2009_2}.}
Since $t\mapsto\int_0^t \langle p, x(s)-z\rangle ds$
    converges in $\R$ (Lemma~\ref{19_1_2009_1}), 
$$
\liminf_{t\to +\infty} \langle  -p,  x(t) -z \rangle\leq 0.
$$
Now take a sequence $(t_n)$, $t_n\to +\infty$, such that
$$
\lim_{n\to +\infty} \langle  -p,  x(t_n) -z \rangle=l
$$
for some $l\in \R$. Since the trajectory $x$ is bounded, and up to a
subsequence,
$$
x(t_n)\rightharpoonup x_{\infty}
$$
for some $ x_{\infty}\in H$. Then
$$
\lim_{n\to +\infty} \langle  p,  x(t_n) -z \rangle= \langle  p,  x_{\infty} -z \rangle.
$$
By weak lower semicontinuity of the function $\Psi$
$$
\Psi( x_{\infty})\leq \liminf_{n\to +\infty}\Psi( x(t_n))=0
$$
thus
$$
 x_{\infty}\in \mbox{argmin} \Psi.
$$
Since $p\in N_{\mbox{argmin} \Psi}(z)$,
$$
 \langle  -p,  x_{\infty} -z \rangle\geq 0.
$$
Thus, every limit point of $ \langle  -p,  x(t) -z \rangle$ is
nonnegative, that is
$$
\liminf_{t\to +\infty} \langle  -p,  x(t) -z \rangle\geq 0.
$$

\noindent Now, for every $t$, since $-p\in \partial \Phi(z)$,
$$
\Phi(x(t))\geq \Phi(z)+ \langle - p,  x(t) -z \rangle,
$$
hence
$$
\liminf_{t\to +\infty}\Phi(x(t))\geq\Phi(z).
$$

\begin{lemma} \label{20_1_2009_4}
Let 
$$
E_2(t)= \Phi(x(t))+\beta(t)\Psi(x(t)).
$$
a) The function $E_2$ is absolutely continuous, and, for a.e. $t$,
$$
\dot{E}_2(t)= -|\dot{x}(t)|^2+\dot{\beta}(t)\Psi(x(t)).
$$
b) The function $E_2$ converges, $\lim_{t\to
  +\infty}\Phi(x(t))=\Phi(z)$, and $\lim_{t\to +\infty} \beta(t)\Psi(x(t))=0$.

\end{lemma}

{\bf Proof of Lemma~\ref{20_1_2009_4}.} See the proof of
Lemma~\ref{20_1_2009_1} for the proof of part a).
Since $\int_0^{+\infty} \beta(t) \Psi(x(t))dt < +\infty$
(Lemma~\ref{19_1_2009_1}) and $ \dot{\beta}\leq k \beta$ (Assumption
$(\mathcal H_2)$), then $\int_0^{+\infty} \dot{\beta}(t)\Psi(x(t))dt < +\infty$,
which implies the convergence of the function  $E_2$. By
Lemma~\ref{19_1_2009_1}, (ii), 
$$
\lim_{t\to+\infty} \int_0^t \Phi(x(s))-\Phi(z) +
\beta(s)\Psi(x(s))ds  \ \mbox{ exists in }\R,
$$
thus 
$$
\lim_{t\to +\infty} E_2(t)=\lim_{t\to
  +\infty}\Phi(x(t))+\beta(t)\Psi(x(t))=\Phi(z).
$$
Hence
$$
\limsup_{t\to +\infty}\Phi(x(t))\leq \lim_{t\to
  +\infty}\Phi(x(t))+\beta(t)\Psi(x(t))\leq \Phi(z).
$$
Since (Lemma~\ref{20_1_2009_2})
$$
\liminf_{t\to +\infty}\Phi(x(t))\geq \Phi(z),
$$
we deduce
$$
\lim_{t\to +\infty} \Phi(x(t))=\Phi(z),
$$
and 
$$
\lim_{t\to +\infty} \beta(t)\Psi(x(t))=0.
$$

In view of Opial's lemma, since we already proved the convergence of $
\|x(t)-z\|$ for every $z\in S=\mbox{argmin} \Phi$, the proof of
Theorem~\ref{Theweak} is finished with the following lemma.

\begin{lemma}\label{20_1_2009_5} If $x(t_n)\rightharpoonup
  x_{\infty}$, then $x_{\infty}\in S=\mbox{argmin}\{ \Phi|\mbox{argmin}\Psi\}$.
\end{lemma}

{\bf Proof of Lemma~\ref{20_1_2009_5}.}
 By weak lower semicontinuity of the functions $\Phi$ and $\Psi$
\begin{eqnarray*}
\Psi( x_{\infty})&\leq&\liminf_{n\to +\infty}\Psi(x(t_n))=0 \ \mbox{ thus
} x_{\infty} \in \mbox{argmin} \Phi.\\
\Phi(x_{\infty})&\leq&\liminf_{n\to +\infty}\Phi(x(t_n))=\Phi(z) \ \mbox{ thus
} x_{\infty} \in \mbox{argmin}\{ \Phi|\mbox{argmin}\Psi\}.
\end{eqnarray*}

\section{Multiscale aspects}
In this section, we show that the system (where $\beta(t)\to +\infty$ as $t\to +\infty$)
$$
(MAG) \qquad \dot{x}(t) +  \partial \Phi (x(t)) +  \beta (t) \partial
\Psi (x(t))\ni 0,
$$
after time rescaling, can be equivalently rewritten  as
$$
(MAG)_{\varepsilon} \qquad \dot{x}(t) +   \partial \Psi (x(t)) +  \varepsilon(t) \partial \Phi (x(t)) \ni 0,
$$
with a positive control $t\mapsto \varepsilon(t)$ that converges to $0$ as
$t \to \infty$. 

By taking $\Phi (x)=\|x\|^2/2$, $(MAG)_{\varepsilon}$ system amounts to the $(SDC)_{\varepsilon}$ system (steepest descent
with control)

$$
(SDC)_{\varepsilon} \qquad \dot{x}(t) +  \partial \Psi (x(t)) +  \varepsilon(t) x(t)   \ni 0.
$$

The $(SDC)_{\varepsilon}$ system plays an important role in optimization and game theory as well as in asymptotic control theory and the study of ill-posed problems. It can be viewed as a Tikhonov-like dynamical system. 
In the particular setting of $(SDC)_{\varepsilon}$, the asymptotic convergence properties of the trajectories depend on whether
$\varepsilon(\cdot)$ is in $L^1(0,\infty)$ or not.  

For the case
$\varepsilon(\cdot)\notin  L^1$, the first general convergence results  go back to \cite{Re}
(based on previous work by \cite{Browd}) and require in addition $\varepsilon(\cdot)$ to be non-increasing. Under these conditions, each trajectory of $(SDC)_{\varepsilon}$ converges strongly to $\bar{x}$, the point of minimal norm in $\mbox{argmin} \Psi$. This case is often referred to as the slow parametrization case (slow convergence of $\varepsilon(\cdot)$ to zero).
In a recent contribution to this subject \cite{CPS}, it is proved that this convergence result still holds without assuming $\varepsilon(\cdot)$ to be non-increasing.

By contrast, for the case
$\varepsilon(\cdot)\in  L^1$, each  trajectory of $(SDC)_{\varepsilon}$ weakly converges to
some point in $\argmin \Psi$ (which depends on the trajectory), a result which is in the line of Bruck theorem. This is ususally referred to as the fast parametrization case.

Among the many papers devoted to $(SDC)_{\varepsilon}$ and related systems, let us   mention 
  \cite{ACom},  \cite{Pey}, \cite{Hir}, \cite{CPS}.
In \cite{ACza} and~\cite{Cza} the authors
  show similar properties concerning  the second order system
$$
(HBFC) \quad \ddot{x}(t)+\gamma \dot{x}(t)+\nabla \Psi(x(t)) + \varepsilon(t) x(t)=0
$$
with $\gamma>0$.  

 Let us now return to the connection with  $(MAG)$. As we shall see, the $(MAG)$ system studied in this paper
leads to an equivalent $(MAG)_{\varepsilon}$ system  with a corresponding control $\varepsilon(\cdot)$ which automatically
satisfies
$$
\varepsilon(\cdot)\notin  L^1.
$$
The strong convergence of the trajectories of $(MAG)$ with
$$
\Phi(x)=\frac{1}{2}\|x\|^2
$$
and hence of $(SDC)_{\varepsilon}$, is a consequence of the strongly monotone case given by
Theorem~\ref{Thsm}. When $\Phi$ is not strongly monotone, the
techniques of the preceding papers remain useless! Indeed, the study of 
$(MAG)$ makes the situation clearer, and, for general $\Phi$ and $\Psi$, allows to find conditions permitting to obtain (weak) convergence to a
point in $S=\mbox{argmin}\{\Phi|\mbox{argmin}\Psi\}$ . This is exactly
condition $(\mathcal H_1)$, which, in terms of the function $\varepsilon(\cdot)$
for the system $(MAG)_{\varepsilon}$,   in classical cases, can be seen as an $L^2$ integrability
condition. For example, if $\Psi(\cdot) =  \frac{1}{2} dist^2 (., C)$, 
the condition will turn to be exactly
$$
\varepsilon(\cdot) \in L^2 \mbox{ and } \varepsilon(\cdot) \notin L^1 .
$$ 
Let us now state precisely the equivalence between the formulations $(MAG)$
and $(MAG)_{\varepsilon}$, and also between $(MAMI)$ and $(MAMI)_{\varepsilon}$.

\begin{lemma}[dictionary]\label{6_2_2009} Let $T_\beta$ and $T_\varepsilon$ be two
  elements in $\left(\R_+\setminus\{0\}\right)\cup{+\infty}$. Take two functions of
  class $C^1$
\begin{eqnarray*}
\beta : [0,T_\beta)&\to &\R_+\setminus\{0\};\\
\varepsilon: [0,T_\varepsilon)&\to &\R_+\setminus\{0\}.
\end{eqnarray*}
Define $ t_\beta:  [0,T_\varepsilon)\to [0,T_\beta)$  \ and \ $
    t_\varepsilon: [0,T_\beta) \to [0,T_\varepsilon)$ by
$$
\int_0^{t_\beta(t)} \beta(s)ds =t \mbox{ and }
\int_0^{t_\varepsilon(t)} \varepsilon(s)ds =t
.
$$
Assume that, for every $t$, 
$$
\varepsilon(t)\beta(t_\beta(t))=1.
$$
Then 
\begin{eqnarray*}
t_\varepsilon \circ t_\beta={\rm id}_{[0,T_\varepsilon)} & and &
  T_\varepsilon=\int_0^{T_\beta} \beta;\\
t_\beta \circ t_\varepsilon={\rm id}_{[0,T_\beta)} & and &
  T_ \beta=\int_0^{T_\varepsilon}\varepsilon;
\end{eqnarray*}
if $x$ is a strong solution of 
$$
(MAG) \qquad \dot{x}(t) +  \partial \Phi (x(t)) +  \beta (t) \partial
\Psi (x(t))\ni 0,
$$
then $x\circ t_\beta$ is a strong solution of 
$$
(MAG)_{\varepsilon} \qquad \dot{w}(t) +  \varepsilon(t) \partial \Phi (w(t)) +            \partial
\Psi (w(t))\ni 0;
$$
conversely, if $w$  is a strong solution of $(MAG)_{\varepsilon}$,
then $w\circ t_\varepsilon$ is  a strong solution of $(MAG)$.

Now, if $x$ is a strong solution of 
$$
 (MAMI) \qquad \dot{x}(t) +  A(x(t)) +  \beta (t) \partial \Psi
 (x(t))\ni 0,
$$
then $x\circ t_\beta$ is a strong solution of 
$$
 (MAMI)_{\varepsilon} \qquad \dot{w}(t) +  \varepsilon(t) A(w(t)) +  \partial \Psi
 (w(t))\ni 0;
$$
conversely, if $w$  is a strong solution of $ (MAMI)_{\varepsilon}$,
then $w\circ t_\varepsilon$ is  a strong solution of $ (MAMI)$.

\end{lemma}

{\bf Proof of Lemma~\ref{6_2_2009}.}  Since $\int_0^{t_\beta(t)} \beta(s)ds =t$, then $t_\beta(\cdot)$
is of class $C^1$ and $\dot{t}_\beta(t)\beta(t_\beta(t))=1$. 
Thus $\dot{t}_\beta(t)=\varepsilon(t)$, and $ t_\beta(t_\varepsilon(t))=
\int_0^{t_\varepsilon(t)} \varepsilon(s)ds =t$. We are now prepared to make the change of variable associated with function $t_\beta(\cdot)$.
Let $x(\cdot)$ be  a strong solution
of $(MAG)$ and write the system
$(MAG)$ at the point ${t_\beta(t)}$:
$$
 \dot{x}(t_\beta(t)) +  \partial \Phi (x(t_\beta(t))) +  \beta (t_\beta(t) ) \partial
\Psi (x(t_\beta(t) ))\ni 0.
$$
After multiplication  by $\dot{t}_\beta(t)$ we get
$$
\dot{t}_\beta(t) \dot{x}(t_\beta(t)) + \dot{t}_\beta(t) \partial \Phi (x(t_\beta(t))) +  \dot{t}_\beta(t)\beta (t_\beta(t) ) \partial
\Psi (x(t_\beta(t) ))\ni 0.
$$
Set $w=x\circ t_\beta$. The map $w=x\circ t_\beta$ is absolutely continuous, and according to $\dot{w}(t)= \dot{t}_\beta(t)\dot{x}(t_\beta(t))$,
$\dot{t}_\beta(t)=\varepsilon(t)$ and $\dot{t}_\beta(t)\beta(t_\beta(t))=1$, we obtain
$$
\dot{w}(t) +  \varepsilon(t) \partial \Phi (w(t)) +            \partial
\Psi (w(t))\ni 0,
$$
that's $(MAG)_{\varepsilon}$.
The notion of strong solution, as given in
Section~\ref{section_sub}, remains valid. Similar arguments work with
$(MAMI)$.
$\Box$

\medskip

Accordingly, all our results can be written for the systems
$(MAG)_{\varepsilon}$ and $ (MAMI)_{\varepsilon}$ . Before doing so, let us analyse the
corresponding assumption $(\mathcal H_1)$.

\begin{remark}\label{r1}{About Assumption ($\mathcal H_1$).} Take $\beta(\cdot)$ and
  $\varepsilon(\cdot)$ as in Lemma~\ref{6_2_2009}, with
  $T_\beta=T_\varepsilon=+\infty$.
  Then, by making the change of variable $t =  t_\beta(s)$ in the following integral, and by using $\dot{t}_\beta(t)\beta(t_\beta(t))=1$, we obtain

\begin{eqnarray*}
\int_{0}^{+\infty} \beta (t) \left(\Psi^*  \left(\frac{p}{ \beta (t)}\right) -
\sigma_C  \left(\frac{p}{ \beta (t)}\right)\right)dt & = &  \int_{0}^{+\infty} \dot{t}_\beta(s)\beta (t_\beta(s)) \left(\Psi^* \left(\frac{p}{ \beta (t_\beta(s))}\right) -
\sigma_C \left(\frac{p}{ \beta (t_\beta(s))}\right)\right)ds\\
&=&\int_{0}^{+\infty} \Psi^*
(\varepsilon(s)p ) - \sigma_C (\varepsilon(s)p)ds.
\end{eqnarray*}

\noindent Thus condition ($\mathcal H_1$) becomes

$$\forall p\in R(N_C) \ \ \int_{0}^{+\infty} \Psi^*
(\varepsilon(t)p ) - \sigma_C (\varepsilon(t)p)dt < + \infty.$$

\noindent In  particular, when  $\Psi (x) =  \frac{1}{2} dist^2 (x, C)$,
then $\Psi^* (x) =  \frac{1}{2} \|x\|^2  +  \sigma_C (x)$ and 
 
$$(\mathcal H_1) \Longleftrightarrow  \int_{0}^{+\infty} \frac
{1}{\beta (t)} dt < + \infty\Longleftrightarrow
\int_{0}^{+\infty}\varepsilon(t)^2 dt < + \infty. \ \Box$$ 
\end{remark}
 
\begin{Theorem} \label{Theweak_eps}
Let  $\Phi$ and $\Psi$ satisfy the assumptions of
Theorem~\ref{Theweak}. Let us assume that,

\begin{itemize}
\item
${(\mathcal H_{1})}_ {\varepsilon}$\mbox{ }\mbox{ }  $ \forall p\in R(N_C),$
  $\displaystyle \int_{0}^{+\infty} \Psi^*
(\varepsilon(t)p ) - \sigma_C (\varepsilon(t)p)dt < + \infty.$

\item
${(\mathcal H_{2})}_ {\varepsilon}$\mbox{ }\mbox{ } $\varepsilon(\cdot)$ is
 a non increasing function of class $C^1$,
  such that $\lim_{t\to+\infty} \varepsilon (t)=0$, $\displaystyle \int_{0}^{+\infty}  \varepsilon (t)dt=+\infty$,
and for some $k\geq 0$, $-k \varepsilon^2\leq 
  \dot{\varepsilon}$.

\end{itemize}

Let $x(\cdot)$ be a strong solution of $(MAG)_{\varepsilon}$. Then:

\begin{eqnarray*}
(i)&\mbox{ weak convergence }&\exists x_{\infty}\in S=\mbox{argmin}\{
  \Phi|\mbox{argmin}\Psi\},\qquad  w-\lim_{t\to+\infty}
  x(t)=x_{\infty};\\
(ii)&\mbox{ minimizing properties } &\lim_{t\to+\infty} \Psi(x(t))=
  0;\\
&&\lim_{t\to +\infty} \Phi(x(t)) = \min\Phi |_{\mbox{argmin}\Psi};\\
(iii)&&\forall z\in S  \lim_{t\to+\infty} \|x(t) -z\| \mbox{  exists  };\\
(iv)&\mbox{ estimations }&\lim_{t\to +\infty} \frac{1}{\varepsilon(t)}\Psi(x(t))=0;\\
&&  \int_{0}^{+\infty} \Psi (x(t)) dt < +\infty;\\
&&\limsup_{\tau\to +\infty} \int_{0}^{\tau}\varepsilon(t)\left(\Phi (x(t))- \min\Phi |_{\mbox{argmin}\Psi}\right)  dt < +\infty.
\end{eqnarray*}
\end{Theorem}

\textbf{{Proof of Theorem} \ref{Theweak_eps}}
Equivalence between Theorem~\ref{Theweak_eps}
and  Theorem~\ref{Theweak} is a 
consequence of Lemma~\ref{6_2_2009}, Remark~\ref{r1} and the equivalent formulation of condition $(\mathcal H_{2})$ as given below:

Write $(\mathcal H_{2})$ at the point $t_{\beta}(t)$ 

$$
\dot{\beta}(t_{\beta}(t))\leq k \beta(t_{\beta}(t))
$$

and multiply by $\dot{t_{\beta}}(t)$ (which is nonnegative)

$$
\dot{\beta}(t_{\beta}(t))\dot{t_{\beta}}(t) \leq k \beta(t_{\beta}(t))\dot{t_{\beta}}(t).
$$
  
Owing to $\beta(t_{\beta}(t))\dot{t_{\beta}}(t) = 1$, we get

$$
\frac{d}{dt} \beta(t_{\beta}(t)) \leq k,
$$

which, by using $\varepsilon(t)\beta(t_\beta(t))=1$, finally yields

$$
\frac{d}{dt} \left(\frac{1}{\varepsilon(t)} \right) =  - \frac{\dot{\varepsilon}(t)}{{\varepsilon}^2(t)}  \leq k
$$

that's  ${(\mathcal H_{2})}_ {\varepsilon}$. $\Box$

\section{Further convergence results, without the growth condition $\dot{\beta}\leq k\beta$}

 It worth noticing that, in the proof of  weak convergence theorem ~\ref{Theweak},
    the growth condition $\dot{\beta}\leq k\beta$  \ is  only used ultimately, in  lemma ~\ref{20_1_2009_4}.
 It allows to develop an energetic argument involving $E_2$. 
 As a byproduct of this technic, one  obtains too that the trajectories have finite kinetic energy: $\int_0^{+\infty} |\dot{x}(t)|^2dt < +\infty.$ 
  It is a natural question to ask whether  theorem ~\ref{Theweak} still holds without  this growth condition. 
  As a positive answer,  when $\mathcal H$ is a finite dimensional space,  we are going to prove that
  convergence of the trajectories holds true without this growth condition.
  Indeed, having in view applications to possibly infinite dimensional problems (PDE's, control), we consider  the more general situation  with functionals $\Phi$ or $\Psi$ which are supposed to be inf-compact. Let us recall that a function $\phi: \mathcal H \rightarrow \R \cup \{+\infty\}$ is said to be  inf-compact if, for every $R>0$ and $l\in \mathbb R$ the lower level set
  
  $$ \left\{ x\in \mathcal H: |x| \leq R , \ \phi(x) \leq l\right\} \ \mbox{ is relatively compact in}    \   \mathcal H. $$

 Our proof is  close to Baillon-Cominetti argument developed in \cite{BC} theorem 2.1. It  relies mainly on topological arguments.
  It turns out that it is more convenient to work with the $(MAG)_{\varepsilon}$ version of our dynamics (see section 4)
   
$$
(MAG)_{\varepsilon} \qquad \dot{x}(t) +   \partial \Psi (x(t)) +  \varepsilon(t) \partial \Phi (x(t)) \ni 0,
$$

with a positive control $t\mapsto \varepsilon(t)$ that converges to $0$ as
$t \to \infty$. Then, the result can be easily  converted in terms of $(MAG)$.

\begin{Theorem} \label{rel-compact-thm}
Let  $\Phi$ and $\Psi$ satisfy the assumptions of
Theorem~\ref{Theweak}. Let us assume that,

\begin{itemize}
\item
${(\mathcal H_{1})}_ {\varepsilon}$\mbox{ }\mbox{ }  $ \forall p\in R(N_C),$
  $\displaystyle \int_{0}^{+\infty} \Psi^*
(\varepsilon(t)p ) - \sigma_C (\varepsilon(t)p)dt < + \infty.$
\item
${(\mathcal H_{2})}_ {\varepsilon}$\mbox{ }\mbox{ } $\varepsilon(\cdot)$ is
 a non increasing function of class $C^1$,
  such that $\lim_{t\to+\infty} \varepsilon (t)=0$, $\displaystyle \int_{0}^{+\infty}  \varepsilon (t)dt=+\infty$.
\item   $(\mathcal H_{3})$ \mbox{ }\mbox{ } $\Phi$ or  $\Psi$ is  inf-compact  and $S$ is bounded.
\end{itemize}

Let $x(\cdot)$ be a trajectory solution of $(MAG)_{\varepsilon}$. Then the following convergence result holds

\begin{equation*}
\exists x_{\infty}\in S=\mbox{argmin}\{
  \Phi|\mbox{argmin}\Psi\},\qquad  \lim_{t\to+\infty}
  x(t)=x_{\infty}
\end{equation*}
the convergence being taken in the strong sense when $\Psi$ is inf-compact and in the weak sense when $\Phi$ is inf-compact.
\end{Theorem}

\textbf{{Proof of Theorem}~\ref{rel-compact-thm}} \ In view of
Lemma~\ref{6_2_2009}, one can easily check that the assumptions of
Theorem~\ref{Theweak} are satisfied, except for the growth condition
$-k \varepsilon^2\leq \dot{\varepsilon} \quad /\quad \dot{\beta}\leq
k\beta$. So all the results in the proof of Theorem~\ref{Theweak}
hold, except for  Lemma~\ref{20_1_2009_4}. 

\textit{Weak convergence of the trajectories.}  We know that, for every $z$ in $S=\mbox{argmin}\{\Phi|\mbox{argmin}\Psi\}$,  
  $\mbox{ } \lim_{t\to+\infty} \|x(t) -z\|$ exists (Lemma~\ref{19_1_2009_1}).
  In order to obtain weak convergence of the trajectory we use Opial lemma. Thus,  we just need to prove  that
  every  weak-limit point of  $x(\cdot)$ belongs to $S$.
  In turn, this will be a straight consequence of 
\begin{equation}\label{relc0} 
   dist(x(t), S) \to 0  \ \mbox{ as } \ t\to + \infty
\end{equation}
 and of the weak lower semicontinuity of the convex continous function  $dist(\cdot, S).$
 In order to prove (\ref{relc0}) let us introduce $h(t) :=  \frac{1}{2}dist(x(t), S)^2$, and  estimate $\dot{h}(t)$. Set 
   $ D_S(x) = \frac{1}{2}dist(x, S)^2$,  which is  convex and differentiable  with $\nabla D_S(x)= x - P_S(x)$, where $P_S$ is the projection of $x$ onto $S$. Then
 for almost every $t$,  
\begin{equation}\label{relc1}
\dot{h}(t) =  \left\langle  \dot{x}(t), x(t) -P_Sx(t)\right\rangle.     
\end{equation}       
Let us rewrite $(MAG)_{\varepsilon}$ as  
 \begin{equation}\label{relc2} 
   \dot{x}(t) +  \varepsilon(t) \partial {\phi}_t (x(t))  \ni 0
\end{equation} 
 with 
 $$ {\phi}_t (x) :=  \Phi (x) +  \frac{1}{\varepsilon(t)}  \Psi (x).  $$
 For almost every $t$,  by using (\ref{relc1}), (\ref{relc2}) and  convexity of $ {\phi}_t$ 
\begin{eqnarray*}
        \dot{h}(t)  & = &  \varepsilon(t) \langle  \partial {\phi}_t (x(t)), P_Sx(t) - x(t)\rangle   \\
                    & \leq &   \varepsilon(t) \left(  {\phi}_t (P_Sx(t)) -  {\phi}_t (x(t))\right). 
\end{eqnarray*}
Let us notice that   ${\phi}_t (P_Sx(t))$ is independent of $t$ (recall that $\Psi =0$  on $S$). It is equal to the optimal value of the limit equilibrium problem, we set 
\begin{center}
$ v_{opt} := \Phi(z)$ \ for all $ z\in  S=\mbox{argmin}\{\Phi|\mbox{argmin}\Psi\}$.
\end{center}
Hence,  
\begin{equation*}\label{relc3} 
   \dot{h}(t) +    \varepsilon(t) \left( {\phi}_t (x(t))  -   v_{opt} \right) \leq 0. 
\end{equation*}
Integrating, and using that the function $h$ is bounded, we get
\begin{equation*}\label{relc4} 
   \int_{0}^{+\infty}  \varepsilon (t) \left( {\phi}_t (x(t))  -   v_{opt} \right) dt < + \infty, 
\end{equation*}
and since $\int_{0}^{+\infty}  \varepsilon (t)dt=+\infty$ we deduce that
\begin{equation*}\label{relc5} 
  \liminf_{t\to +\infty} \left( {\phi}_t (x(t))  -   v_{opt} \right) \leq 0. 
\end{equation*}
By definition of $\liminf$, this implies the existence of a sequence $t_k \to +\infty$ such that $\lim_{k\to +\infty}  {\phi}_{t_k} (x(t_k))  \leq   v_{opt} .$
From  boundedness of the trajectory $x(\cdot)$ and  inf-compactness assumption $(\mathcal H_{3}),$ we deduce that the sequence $(x(t_k))$ is relatively compact in $\mathcal H.$
On the other hand, as $t\to + \infty$, the sequence of functions ${\phi}_t$ converges increasingly  to $ \Phi + \delta_C $ where $C= \mbox{argmin}\Psi$. By a classical result, monotone convergence implies  epiconvergence ($\Gamma$-convergence) (see \cite{A} theorem 2.40) with its accompanying variational properties.  Still denoting $x(t_k)$ a subsequence  which converges to some $x_{\infty}$,
 we obtain that $x_{\infty} \in \mbox{argmin}\{\Phi + \delta_C\} = S$. From
\begin{equation*}\label{relc6} 
  h(t_k) = \frac{1}{2}dist(x(t_k), S)^2     \leq \frac{1}{2} |x(t_k) - x_{\infty}|^2, 
\end{equation*}
we obtain the convergence of $ h(t_k)$ to zero.
Thus, we have established the existence of a sequence $t_k \to +\infty$ such that $ h(t_k)$ tends to zero. The proof will be completed by proving that $ \lim_{t\to +\infty}h(t)$ exists.
Let us return to (\ref{relc3}), introduce the set  
$$
  S_t := \{x: {\phi}_t (x) \geq v_{opt}\}
$$ 
 and observe that 
\begin{equation*}\label{relc7} 
  x(t) \in S_t \Rightarrow \dot{h}(t) \leq 0.
\end{equation*}
 Set 
\begin{equation}\label{relc71} 
b(t):= \mbox{sup} \{D_S(x):  {\phi}_t (x) \leq v_{opt} \}.
\end{equation}
  Thus, 
 \begin{equation*}\label{relc8} 
  D_S (x) > b(t)   \Rightarrow x\in  S_t                                                      
\end{equation*}
 and
 \begin{equation*}\label{relc9} 
 h(t)= D_S (x(t)) > b(t)   \Rightarrow x(t) \in  S_t    \Rightarrow \dot{h}(t) \leq 0.                                                  
\end{equation*}
Let us show that $\mbox{lim}_{t\to +\infty} b(t) =0$, with $b(t)$ defined in (\ref{relc71}).
Let us argue by contradiction and suppose that there exists $\varepsilon > 0$, $t_k \to +\infty$, and $x_k \in \mathcal H$ with ${\phi}_{t_k} (x_k) \leq v_{opt}$
and $dist (x_k, S) > \varepsilon$. By convex combination of  $x_k$ with $P_S (x_k)$ we can assume, without loss of generality, that $dist (x_k, S) = \varepsilon$.
The set $S$ has been assumed to be bounded. It follows that the sequence  $(x_k)$ is bounded.  From the inf-compactness assumption $(\mathcal H_{3})$ and ${\phi}_{t_k} (x_k) \leq v_{opt}$
we deduce that the sequence $(x_k)$ is relatively compact. Thus, passing to a subsequence, we can assume that $x_k$ strongly converges to some $\bar{x}$. From ${\phi}_{t_k} (x_k) \leq v_{opt}$
and the preceding epi-convergent argument we  get $\bar{x} \in S$. This clearly contradicts the fact that $dist (x_k, S) = \varepsilon$ and $x_k$ strongly converges to  $\bar{x}$.\\
Collecting the preceding results, we are in position to apply the
following lemma   to obtain that $ \lim_{t\to +\infty}h(t)$ exists, and
hence conclude to the weak convergence.  We give the proof of the lemma for the convenience of the reader.

\begin{lemma}\label{Bai-Com} (\cite{BC} lemma 2.2) Let  $h(\cdot)$ \ and $b(\cdot)$ be two real-valued functions defined on $(0, +\infty)$ with $h(\cdot)$  absolutely continuous and nonnegative. Let us assume

\begin{itemize}
\item (i)  \ $b(t) \to 0  \ \mbox{ as } \ t\to  +\infty$;
\item(ii) there exist a set $N \subset (0, +\infty)$ of zero Lebesgue measure such that
$h(t) > b(t) \Rightarrow  \dot{h}(t)  \leq 0 \mbox{ for all t}\notin  N$.
\end{itemize}

Then h(t) has a limit for  $t\to  +\infty.$
\end{lemma}

{\bf Proof of Lemma~\ref{Bai-Com}.} Replacing $b(t)$ by sup$\left\{b(s): s\geq t\right\}$, we can assume, without loss of generality, that $b(\cdot)$ is nonincreasing.
Since by assumption $h(\cdot)$ is nonincreasing when above $b(\cdot)$, it follows easily that, if $h(s) \leq b(s)$ for a given $s\geq 0$, then $h(t) \leq b(s)$ for all $t \geq s$.

Let us examine the two situations:
Either  there exists some $\bar{s}$ with $h(t) > b(t)$ for all $t\geq \bar{s}$, then  $h(\cdot)$ is nonincreasing over $\left[\bar{s}, +\infty\right)$ and therefore it converges when $t\to +\infty$.
Or there exists a sequence $s_k \to +\infty$ such that $h(s_k) \leq b(s_k)$ for all $k\in \mathbb N$. By using the above remark, we deduce  $0 \leq h(t) \leq b(s_k)$ for all $t\geq s_k$.
Since $b(s_k) \to 0$, we conclude that $h(t) \to 0$ as $t\to +\infty.$ $\Box$
  
\textit{Strong convergence when $\Psi$ is inf-compact.}  Since
the trajectory $x(\cdot)$ is bounded, and since $\Psi(x(t))\to 0$, we
conclude that the trajectory $x(\cdot)$ is relatively compact. It weakly
converges, thus strongly converges. $\Box$

\section{Applications}

We first show how our study fits coupled gradient dynamics; then, we consider two particular situations, firstly domain decomposition for elliptic PDE's, secondly best response dynamic approach to Nash equilibria for potential games.

\subsection{Coupled gradient dynamics}

Throughout this section we make the following assumptions: 
\begin{itemize}
\item $\mathcal H = \Xiu \times \Xid$ is the cartesian product of two Hilbert spaces, set $x=(x_1,x_2)$;
\item  $\Phi(x) = f_1(x_1) + f_2(x_2) + \phi (x_1,x_2) $,  \  $f_1\in \Gamma_0 (\Xiu)$,  $f_2\in \Gamma_0 (\Xid)$ 
are closed convex functions, $\phi : \Xiu \times \Xid \to \R$ is a smooth convex coupling function;
\item $\Psi (x) = \frac{1}{2}\|L_1x_1 - L_2x_2\|^2_{\Z}$,  \   
$L_1 \in L(\Xiu, \Z)$ and  $L_2 \in L(\Xid, \Z)$ are linear continuous operators acting respectively from $\Xiu$  and $\Xid$ into a third Hilbert space $\Z$;
\item $\beta : \R^+ \rightarrow \R^+$ is  a function of $t$ which tends to $+\infty$ as $t$ goes to $+\infty$.
\end{itemize}

\noindent In this setting, $(MAG)$ system
$$ \qquad \dot{x}(t) +  \partial \Phi (x(t)) +  \beta (t) \partial \Psi (x(t))\ni 0$$

\noindent becomes

\begin{equation}\label{couplsys}
\left\{
\begin{array}{l}
\dot{x_1}(t) +  \partial f_1(x_1(t)) + {\nabla}_{x_1}\phi (x_1 (t),x_2 (t))  +  \beta (t)L_1^{*} ( L_1x_1(t) - L_2x_2(t))\ni 0\\
\rule{0pt}{25pt}
\dot{x_2}(t) +  \partial f_2(x_2(t)) + {\nabla}_{x_2}\phi (x_1 (t),x_2 (t)) + \beta (t) L_2^{*} ( L_2x_2(t) - L_1x_1(t))\ni 0.
\end{array}\right.
\end{equation}
 
\noindent Because of the quadratic property of $\Psi$, condition $({\mathcal H}_1)$ can be equivalently written 

$$ \int_{0}^{+\infty} \frac {1}{\beta (t)} dt < + \infty.$$

\noindent As a straight application of theorem \ref{Theweak}, assuming $({\mathcal H}_1)$ and the growth condition 

$$\dot{\beta}\leq k\beta$$

\noindent we obtain that
$x(t)= (x_1(t),x_2(t))\rightarrow x_{\infty}=(x_{1,\infty},x_{2,\infty})$ weakly in  $\mathcal H$
where $(x_{1,\infty},x_{2,\infty})$ is a solution of

\begin{equation}\label{structopt}
\mbox{min} \left\{ f_1(x_1) + f_2(x_2) + \phi (x_1,x_2): \mbox{ }\mbox{ }L_1x_1 - L_2x_2 = 0\right\}.
\end{equation}

\noindent In theorems \ref{Thsm} and \ref{rel-compact-thm} we describe several situations where similar conclusions hold without condition $\dot{\beta}\leq k\beta$.

\noindent Structured optimization problems  (\ref{structopt}) occur in various domains:

\begin{itemize}

\item In game theory (see subsection 6.3 for further details), (\ref{structopt}) describes Nash equilibria of the potential game (here team game) with two players $1$,  $2$ and respective
static loss functions: 
$$\left\{
\begin{array}{l}
\vspace{0.2cm}
F_1:(x_1,x_2)\in \Xiu \times \Xid \rightarrow F_1(x_1,x_2) = f_1(x_1) +  \phi(x_1,x_2) \mbox{ } \mbox{ if }\mbox{ }  L_1x_1 - L_2x_2 = 0 , \mbox{ }+\infty  \mbox{ elsewhere}  \\ 
\rule{0pt}{15pt}
F_2:(x_1,x_2)\in \Xiu \times \Xid \rightarrow F_2(x_1,x_2) = f_2(x_2) +  \phi(x_1,x_2) \mbox{ }\mbox{ if } \mbox{ } L_1x_1 - L_2x_2 = 0, \mbox{ }+\infty \mbox{ elsewhere}. 
\end{array}\right.$$

\noindent The $f_i (.)$ represent the individual payoffs of the players, $\phi(.,.)$ is their joint payoff, and $L_1x_1 - L_2x_2 = 0$ is a  constraint expressing, in a normalized form,   global limitation of the resources.
In that case, a discrete version of system (\ref{couplsys}) provides  a best response dynamic approach to such equilibria.

\vspace{0.2cm}

\item In optimal control theory, the constraint $\mbox{ }\left\{L_1x_1 = L_2x_2\right\}\mbox{ }$ is the state equation, which relates the state variable $\mbox{ }$$x_1$$\mbox{ }$
to the corresponding control  $\mbox{ }$$x_2$$\mbox{ }$. In that case, the criteria which is to minimize naturally splits into the sum of two costs:  the cost to be far from a desired state  and the cost of the control.

\vspace{0.2cm}

\item Variational formulation of phase transition, cracks and fissures, image segmentation and many others naturally lead to minimization problems of type (\ref{structopt}). Here, functions $f_i$ represent the internal energy 
of the different phases, while  coupling terms represent  energies located at the interfaces, as well as some  transmission conditions. Indeed, for numerical purpose, it is also interesting  to relax the classical variational formulation of elliptic boundary value problems into the form (\ref{structopt}) in order to perform decomposition of domains methods (see section 6.1 below).
 
\end{itemize}

\vspace{0.1cm} 

\noindent An important question is the modelling of  dynamic approach to such equilibria, and the design of iterative numerical schemes (algorithms) for solving the corresponding problems.  Both concern   asymptotic behavior of associated discrete dynamical systems. Indeed, as a general rule, continuous versions of these systems offer  flexible  tools allowing  a deeper understanding of their mathematical properties. Moreover, they  may suggest extensions and connections with other domains. Let us illustrate this with two applications in different domains: one concerns domain decomposition methods for PDE's, and the other best response dynamics in potential games.
\subsection{Domain decomposition for PDE's}
In the following example, the domain $\Omega$ naturally splits into two elementary non overlapping subdomains  ${\Omega}_i$  (i=1,2) with common interface $\Gamma,$ i.e., 
$\Omega = {\Omega}_1  \cup {\Omega}_2  \cup  \Gamma $.
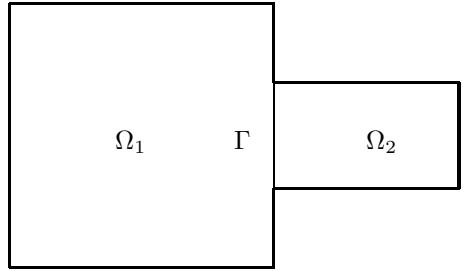
\vspace{0.3cm}

\begin{center}
\begin{picture}(200,100)(0,0)
\thicklines
\put(0,0){\line(1,0){100}}
\put(100,0){\line(0,1){30}}
\put(100,30){\line(1,0){70}}
\put(170,30){\line(0,1){40}}
\put(170,70){\line(-1,0){70}}
\put(100,70){\line(0,1){30}}
\put(100,100){\line(-1,0){100}}
\put(0,100){\line(0,-1){100}}
\thinlines
\put(100,30){\line(0,1){40}}
\put(40,45){$\Omega_1$}
\put(135,45){$\Omega_2$}
\put(85,45){$\Gamma$}
\end{picture}
\end{center}

\noindent  Given $h\in L^{2}(\Omega)$,  Dirichlet problem on $\Omega$  \ consists finding  $u: \Omega \rightarrow \R$ solution of

$$\left\{
\begin{array}{l}
-\Delta u = h \mbox{ }\mbox{ }\mbox{on} \mbox{ }\Omega\\
\rule{0pt}{15pt}
u = 0 \mbox{ }\mbox{ }\mbox{on} \mbox{ }\partial\Omega.\\
\end{array}\right.$$

\noindent Its variational formulation can be equivalently formulated as

$$\min\left\{
  \frac{1}{2}\int_{\Omega_1}|\nabla v_1|^2  - \int_{\Omega_1}hv_1
  +\frac{1}{2}\int_{\Omega_2}|\nabla v_{2}|^2 - \int_{\Omega_2}hv_2
  :\mbox{ }\mbox{ }v = 0 \mbox{ }\mbox{ }\mbox{on} \mbox{ }\partial\Omega, \mbox{ }\mbox{ }\left[v\right] = 0 \mbox{ }\mbox{ }\mbox{ on } \Gamma \right\}
  $$

\noindent where $v = v_i$ on   ${\Omega}_i$ and       $\left[v\right]$ is the  jump of $v$ through the interface $\Gamma.$  
 
\noindent  Indeed, the above problem falls into the setting of (\ref{structopt}) (with $\phi =0$)
  
$$\min\left\{ f_1(v_1)+ f_2(v_2):  v_1\in \Xiu, v_2\in \Xid, \mbox{ }\mbox{ } L_1(v_{1}) - L_2(v_{2}) = 0 \right\}$$

\noindent  by taking

\begin{itemize}
\item $\mathcal X_i = \{v\in\ H^1(\Omega_i): \mbox{ }\mbox{ } v = 0\mbox{ on }\partial\Omega\cap\partial{\Omega_i}\}, \mbox{ }\mbox{ }v = v_{i} \mbox{ }\mbox{on} \mbox{ }\Omega_i, \mbox{ }i=1,2.$

\vspace{0.1cm} 

\item $f_i(v_i) =  \frac{1}{2}\int_{\Omega_i}|\nabla v_i|^2   - \int_{\Omega_i}hv_i.$

\vspace{0.1cm} 

\item $L_i: H^1(\Omega_i)\rightarrow\mathcal Z = L^2(\Gamma) \mbox{ }\mbox{ the Sobolev trace operator}, \ i=1,2.$

\vspace{0.1cm} 

\item  $\left[v\right]$ = $L_1(v_{1}) - L_2(v_{2})$ = jump of $v$ through the interface $\Gamma.$
\end{itemize}

\noindent Let us equip $\mathcal X_i = \{v\in\ H^1(\Omega_i): \  v = 0 \mbox{ on } \partial\Omega\cap\partial{\Omega_i}\}$  \ with the scalar product
 $\left\langle u, v \right\rangle =  \int_{\Omega_i} \nabla u . \nabla v $. By Poincar\'e inequality, the induced norm is equivalent to the usual norm of $H^1(\Omega_i)$. 
 Then, (\ref{couplsys}) reads as follows
 
$$\left\{
\begin{array}{l}
-\Delta \frac{\partial u_1}{\partial t}  - \Delta u_1 = h_1  \mbox{ }\mbox{ }\mbox{on} \mbox{ }\Omega_1\\
\rule{0pt}{30pt}
-\Delta \frac{\partial u_2}{\partial t}  - \Delta u_2 = h_2  \mbox{ }\mbox{ }\mbox{on} \mbox{ }\Omega_2\\
\rule{0pt}{30pt}
\frac{\partial \dot{u_1}(t) }{\partial \nu_1} + \frac{\partial u_1}{\partial \nu_1}(t) - \beta(t) \left[u(t)\right] = 0\mbox{ }\mbox{ } \mbox{on} \mbox{ }\Gamma\\
\rule{0pt}{30pt}
\frac{\partial \dot{u_2}(t) }{\partial \nu_2} + \frac{\partial u_2}{\partial \nu_2}(t) + \beta(t) \left[u(t)\right] = 0 \mbox{ }\mbox{ }\mbox{on} \mbox{ }\Gamma\\
\end{array}\right.$$

\vspace{0.2cm} 

\noindent A standard implicit discretization scheme leads to the following  alternating  algorithm with Dirichlet-Neumann transmission conditions:

\vspace{0.2cm} 

\begin{center}
 $(u_{1,k}, u_{2,k})\rightarrow (u_{1,k+1}, u_{2,k})\rightarrow (u_{1,k+1}, u_{2,k+1})$  \ with  $\beta_k \rightarrow + \infty $;
\end{center}

$$\left\{
\begin{array}{l}
-(1+\alpha) \Delta u_{1,k+1} = h_1 - \alpha\Delta u_{1,k} \mbox{ }\mbox{ }\mbox{on} \mbox{ }\Omega_1\mbox{ }\mbox{ }\\
 \rule{0pt}{15pt}
(1+\alpha) \frac{\partial u_{1,k+1}}{\partial \nu_1} + \beta_k u_{1,k+1} = \beta_k u_{2,k} + \alpha  \frac{\partial u_{1,k}}{\partial \nu_1}\mbox{ }\mbox{ }\mbox{on} \mbox{ }\Gamma\mbox{ }\mbox{ }\\
\rule{0pt}{15pt}
u_{1,k+1} = 0  \mbox{ }\mbox{ }\mbox{on} \mbox{ }\partial{\Omega_1}\cap\partial \Omega\mbox{ }\mbox{ }\\
\end{array}\right.$$
$$\left\{
\begin{array}{l}
\vspace{0.2cm}
- (1+ \alpha)\Delta u_{2,k+1} = h_2 - \alpha \Delta u_{2,k} \mbox{ }\mbox{ }\mbox{on} \mbox{ }\Omega_2\\
\rule{0pt}{15pt}
(1+ \alpha) \frac{\partial u_{2,k+1}}{\partial \nu_2} + \beta_k u_{2,k+1} = \beta_k u_{1,k+1} + \alpha   \frac{\partial u_{2,k}}{\partial \nu_2}\mbox{ }\mbox{ }\mbox{on} \mbox{ }\Gamma\\
\rule{0pt}{15pt}
u_{2,k+1} = 0  \mbox{ }\mbox{ }\mbox{on} \mbox{ }\partial{\Omega_2}\cap\partial \Omega\\

\end{array}\right.$$

In the above algorithm, one has to solve  boundary value problems alternatively on  $\Omega_1$ and $\Omega_2$. Thus the initial problem has been decomposed into more elementary subproblems.
This approach can be advantageously combined with Lagrangian technics,  parallel computing methods and fits well constraints on the data as well as unilateral transmission conditions, as long as convexity properties are satisfied, see  \cite{ASouey} and reference herein.

\subsection{Best response dynamics for potential games}

\vspace{0.2cm}

Consider the potential game (here team game) with two players $1$ and $2$ whose respective
static loss functions are given by 

$$\left\{
\begin{array}{l}
\vspace{0.2cm}
F_1:(x_1,x_2)\in \Xiu \times \Xid \rightarrow F_1(\xi,x_2) = f_1(x_1) +  \phi(x_1,x_2) \mbox{ } \mbox{ if }\mbox{ }\mbox{ }  L_1x_1 - L_2x_2 = 0 , \mbox{ }+\infty  \mbox{ elsewhere}  \\ 
\rule{0pt}{15pt}
F_2:(x_1,x_2)\in \Xiu \times \Xid \rightarrow F_2(x_1,x_2) = f_2(x_2) +  \phi(x_1,x_2) \mbox{ }\mbox{ if } \mbox{ } L_1x_1 - L_2x_2 = 0, \mbox{ }+\infty \mbox{ elsewhere}. 
\end{array}\right.$$
\vspace{0.3cm}

\noindent Because of the particular structure (the joint payoff $\phi(\cdot,\cdot)$ of the two players is the same), the Nash equilibria are the solutions of the convex constrained minimization problem

$$ \mbox{min} \left\{f_1(x_1) + f_2(x_2) + \phi(x_1,x_2):  L_1x_1 - L_2x_2 = 0\right\}. $$

\vspace{0.1cm}

\noindent The constraint $ L_1x_1 - L_2x_2 = 0$ reflects some limitation on the global resources of the agents.

\noindent A central question in game theory, decision sciences and economics is to describe realistic dynamics which converge to such equilibria.
In this context, the corresponding (MAG) dynamic 

$$\left\{
\begin{array}{l}
\dot{x_1}(t) +  \partial f_1(x_1(t)) + {\nabla}_{x_1}\phi (x_1 (t),x_2 (t))  +  \beta (t)L_1^{*} ( L_1x_1(t) - L_2x_2(t))\ni 0\\
\rule{0pt}{15pt}
\dot{x_2}(t) +  \partial f_2(x_2(t)) + {\nabla}_{x_2}\phi (x_1 (t),x_2 (t)) + \beta (t) L_2^{*} ( L_2x_2(t) - L_1x_1(t))\ni 0
\end{array}\right.$$

\noindent provides a  valuable guideline.
Indeed,  discretization of this continuous dynamics leads to the following
``Best reply dynamic with cost to change'', (players $1$ and $2$ play alternatively)

\begin{center}
 $(x_{1,k}, x_{2,k})\rightarrow (x_{1,k+1}, x_{2,k})\rightarrow (x_{1,k+1}, x_{2,k+1})$ with  $\beta_k \rightarrow + \infty $;
\end{center}

 $$ \left\{                                                             
  \begin{array}{l} 
 x_{1,k+1} = \mbox{argmin}
    \{f_1(\xi)+ \phi(\xi,x_{2,k}) + \frac{\beta_k}{2}\|L_1\xi - L_2x_{2,k}\|^2 + \frac{\alpha}{2} \parallel \xi - x_{1,k}\parallel^{2}_{\Xiu }:\mbox{ }
      \xi\in\Xiu \}\\
\rule{0pt}{30pt} 
 x_{2,k+1}=\mbox{argmin}
    \{f_2(\eta)+  \phi( x_{1,k+1},\eta) +  \frac{\beta_k}{2} \|L_1 x_{1,k+1} - L_2\eta\|^2 + \frac{\nu}{2} \parallel\eta - x_{2,k}\parallel^{2}_{\Xid}: \mbox{ }
      \eta\in\Xid\}.\\
   \end{array}\right.$$
   
\vspace{0.2cm}

\noindent For the cognitive and psychological interpretation of the costs to move terms $\parallel \xi - x_{1,k}\parallel^{2}_{\Xiu }$ and
$ \parallel\eta - x_{2,k}\parallel^{2}_{\Xid}$
 consult \cite{ARS}, \cite{ABRS}, \cite{AS}.
 The parameter $\beta_k$ traducts some adaptive  behavior of the agents,
 with endogenous and/or exogenous aspects.
This discrete dynamic provides an elementary model for "how to learn sharing limited resources".

\end{document}